\documentclass[11pt]{article}
\usepackage[letterpaper,textheight=9in,textwidth=6in]{geometry}

\usepackage{amsmath, amsfonts, amssymb, amstext, amsthm} 
\usepackage[affil-it]{authblk} 
\usepackage[titletoc,title]{appendix}
\usepackage{enumitem}
\setlist[enumerate]{leftmargin=.5in}

\usepackage{graphicx,array}
\def\R{\mathbb{R}}
\def\N{\mathbb{N}}

\def\Z{\mathbb{Z}}

\def\Xint#1{\mathchoice
{\XXint\displaystyle\textstyle{#1}}%
{\XXint\textstyle\scriptstyle{#1}}%
{\XXint\scriptstyle\scriptscriptstyle{#1}}%
{\XXint\scriptscriptstyle\scriptscriptstyle{#1}}%
\!\int}
\def\XXint#1#2#3{{\setbox0=\hbox{$#1{#2#3}{\int}$ }
\vcenter{\hbox{$#2#3$ }}\kern-.6\wd0}}

\def\dashint{\Xint-}

\renewcommand{\Re}{\mathrm{Re}\,}

\def\rme{\mathrm{e}}

\def\rmO{\mathrm{O}}

\def\rmd{\mathrm{d}}

\def\rmi{\mathrm{i}}
\newenvironment{Proof}[1][.]%
 {\begin{trivlist}\item[]\textbf{Proof#1 }}%
 {\hspace*{\fill}$\rule{0.4\baselineskip}{0.4\baselineskip}$\end{trivlist}}

\renewcommand{\qed}{\hspace*{\fill}$\rule{.4\baselineskip}{.4\baselineskip}$}

\newtheorem{Lemma}{Lemma}[section]
\newtheorem{Theorem}{Theorem}
\newtheorem{Proposition}[Lemma]{Proposition}
\newtheorem{Corollary}[Lemma]{Corollary}
\newtheorem{Remark}[Lemma]{Remark}

\numberwithin{equation}{section}

\title{Advection-diffusion dynamics with nonlinear boundary flux as a model for crystal growth}
\author{Antoine Pauthier and Arnd Scheel
}

\newcommand{\Addresses}{{
  \bigskip
  \footnotesize

  A.~Pauthier (Corresponding author) and A.~Scheel,\newline \textsc{School of Mathematics, University of Minnesota, 206 Church St. SE, Minneapolis, MN 55455}\par\nopagebreak
  \textit{E-mail address}, A.~Pauthier: \texttt{apauthie@umn.edu}, A.~Scheel: \texttt{scheel@umn.edu}

}}

\begin{document}
 \maketitle
\begin{abstract}
We analyze the effect of nonlinear boundary conditions on an advection-diffusion equation on the half-line. 
Our model is inspired by models for crystal growth where diffusion models diffusive relaxation of a displacement field, 
advection is induced by apical growth, and boundary conditions incorporate non-adiabatic effects on displacement at the boundary. 
The equation, in particular the boundary fluxes,  possesses a discrete gauge symmetry, and we study the role of simple, 
entire solutions, here periodic, homoclinic, or heteroclinic  relative to this gauge symmetry, in the global dynamics. 
\end{abstract}
 

\section{Introduction}
We present results on entire solutions and asymptotic behavior in a very simple model for the effect of boundaries and growth on crystalline phases. 
We first present some background for our model and then describe our main results. 

\paragraph{Derivation of phase-diffusion, strain-displacement, and effective boundaries.}
We are interested in the effect of boundaries and growth on crystalline phases. In the simplest context, such systems can be modeled by a displacement 
field $u$ obeying a gradient flow to an elastic energy $\int|\nabla u|^2$, resulting in a heat equation in a domain $\Omega$. Again idealizing the situation, 
we consider one-dimensional domains $\Omega\subset \R$, and focus on a single boundary, only, $\Omega=\R^+$. Such problems can be derived at various levels of 
rigor in mesoscopic systems for striped phases and are then known as phase-diffusion equations. A prototypical example are modulations of 
periodic patterns in the Swift-Hohenberg equation 
\[
 \partial_t U = -(\partial_{xx}+1)^2 U + \mu U - U^3. 
\]
For $\mu>0$, the equation possesses a family of periodic solutions $U(kx;k)$,  $U(\xi;k)=U(\xi+2\pi;k)$, $k\sim 1$. 
Posing the equation on $x\in\R$, one can approximate solutions to the Swift-Hohenberg equation by modulations of periodic solutions $U(kx+\Phi;k+\Phi_x)$, $\Phi=\Phi(t,x)$, 
provided $\Phi$ solves a phase-diffusion equation 
\[
\partial_t\Phi=d(k)\partial_{xx}\Phi,
\]
where $d(k)$ is an effective diffusivity. Approximation results of this type are usually obtained in a long-wavelength scaling regime, 
$\Phi=\varepsilon \phi(\varepsilon^2 t,\varepsilon x)$, but often seem to hold in much more generality \cite{dsss,ercolani,Noble_Rodrigues_ARMA1,Noble_Rodrigues_inventiones,schneidermelbourne,sssu,schneidernws}. 
We note here that the two-parameter family of equilibria $U(kx+\varphi;k)$, $k\sim 1, \varphi\in [0,2\pi)$ corresponds to the trivial equilibria 
$\Phi(t,x)\equiv \kappa x+\varphi$ in the phase-diffusion equation. Also, note that $\Phi$, as a phase, takes values in the circle $\R/ 2\pi\Z$, 
that is, our phase-diffusion equation possesses a gauge symmetry $\Phi\mapsto \Phi+2\pi$. 

The effect of boundaries on striped phases can be quite complex. In this simple one-dimensional caricature, earlier work demonstrated that boundaries 
act through an effective strain-displacement relation \cite{beekie,ms,weinburd}, relating strain $k$ and displacement $\varphi$. 
In the simplest scenario, the boundary acts as a selector of wavenumbers, $k=g(\varphi)$, with $g(\varphi)=g(\varphi+2\pi)$, 
thus providing effective boundary conditions for the phase-diffusion equation $
\Phi_x=g(\Phi)
$ and selecting stationary solutions $\Phi(x)=g(\varphi)x+\varphi$. 
While we are not aware of rigorous approximation results for temporal dynamics including such boundary conditions, 
they have been found to well predict the interaction of striped phases with boundaries in many scenarios \cite{beekie,avery}.

\paragraph{Apical growth.}
Of particular interest to us are scenarios where the domain occupied by the crystalline phase grows in time. 
At constant growth rate $c$, such a domain would be $\Omega_t=\{x>-ct\}$, which leads, in a comoving frame $\xi=x+ct$ to the equation 
\[
 \partial_t\Phi=\partial_{\xi\xi}\Phi - c\Phi_\xi,\qquad \xi>0.  
\]
In the following we study this equation with nonlinear boundary conditions subject to the gauge symmetry $\Phi\mapsto \Phi+2\pi$. 
Given the simplicity of the equation, we hope that the results are of interest possibly even outside of the context of crystal 
growth and we will not stress the relationship in most of our exposition, only return to this perspective in the discussion.

\paragraph{Advection-diffusion and functional-analytic setup.}
Consider the advection-diffusion equation with non-linear flux,
\begin{equation}
 \begin{cases}\label{pb1}
  \partial_t u = \partial_{xx}u-c\partial_x u, \qquad & x>0,\ t>0\\
  \partial_x u = g(u), \qquad & x=0,\ t>0,
 \end{cases}
\end{equation}
where $c\geq 0$ and  the flux is assumed to be smooth, $g\in C^2\left(\R^,\R\right)$, and to possess the gauge invariance 
\begin{equation}\label{eq:condong}
g(u+2\pi)=g(u).
\end{equation}
We are interested in the long-time behavior of solutions of (\ref{pb1}), in particular $\omega$-limit sets of trajectories, which consist of entire solutions, that is,  solutions defined 
for all $t\in\R$. 
Throughout, we think of solutions as classical solutions and  assume some growth condition.  Namely, denoting $\displaystyle X_+:=BUC\left(\R^+\right)$ 
the space of bounded and uniformly continuous functions on $[0,\infty)$,
we consider the solutions of (\ref{pb1}) in the function space
\begin{equation}\label{def:Yc}
 Y_\mathrm{c}:=\left\{v\ : \ x\mapsto \rme^{-\frac{c}{2}x}v(x)\in X_+\right\},
\end{equation}
and we are investigating solutions of (\ref{pb1}) such that $u(\cdot,t)$ belongs to $Y_\mathrm{c}$, for all $t$. 

In the following, we distinguish positive and zero advection speeds,  $c>0$ and $c=0$. In the case of positive speed, we first discuss the case $g>0$, before turning to the case where $g$ has zeros. 


\paragraph{Positive flux and periodic solutions.}
Suppose  $g>0$,  the case $g<0$ being obtained by reflection $u\mapsto -u$.  Scaling and shifting $u$, we may consider only fluxes 
\begin{equation}\label{def:gtheta}
 g(u)=g(u;\vartheta)=1+\vartheta h(u),\ \textrm{ with }h(\pi)=\min h=-1,\ \max h=1,\ \vartheta\in(-1,1).\tag{$\mathcal{BC}_\vartheta$}
\end{equation}
A prototypical example for such a nonlinear term is
\begin{equation}\label{eq:geqcos}
 g(u;\vartheta) = 1+\vartheta\cos(u).
\end{equation}
One quickly finds that there are no stationary solutions to (\ref{pb1}) satisfying (\ref{def:Yc}). Moreover, $|\vartheta|<1$ and thus $\partial_xu>00$ at $x=0$ suggests 
that solutions decrease monotonically. The expected asymptotic behavior is a time-decreasing solution, which eventually would resonate 
with the gauge symmetry and lead to asymptotic behavior $u(t,x)=u(t+T,x)+2\pi$ for some minimal period $T>0$. 

\begin{Theorem}[Periodic solutions up to gauge]\label{thm1}
 Consider (\ref{pb1}) with boundary condition  (\ref{def:gtheta}) with $c>0$. Then for all $\vartheta\in(-1,1)$, there exists an entire solution $u=u(x,t;\vartheta)$ of (\ref{pb1}) such that
$u$ is 
\begin{enumerate}
 \item relative periodic in time, 
\begin{equation*}
  u(x,t+T;\vartheta)=u(x,t;\vartheta)-2\pi, \textrm{ for some minimal }T=T(\vartheta),
\end{equation*}
\item monotone in space and time, $\partial_x u>0$ and $\partial_t u<0$,
\item unique (up to time translation) in the class of entire solutions with $\partial_t u<0$. 
\end{enumerate}
\end{Theorem}

\paragraph{Sign-changing flux and heteroclinic solutions.}
The condition $|\vartheta|<1$, or equivalently $g>0$, is essential in the proof. In case $|\vartheta|\geq1$, 
(\ref{pb1}) admits constant stationary solutions preventing the existence of such a time-decreasing dynamics. Our first result is concerned with the fate of periodic solutions as $|\vartheta|\nearrow1$.

 \begin{Theorem}[Saddle-node on a limit cycle at $\vartheta=1$]\label{prop2}
 Consider (\ref{pb1}) with boundary condition  (\ref{def:gtheta}) with $c>0$. Suppose that $h(u)>-1$ except for $u=\pi \mod 2\pi$. 
 Then there exists an entire solution $u_1(x,t)$ of (\ref{pb1})--$(\mathcal{BC}_1)$ that is
\begin{enumerate}
 \item homoclinic up to the gauge symmetry, 
\[ u_1(\cdot,t)\underset{t\to-\infty}{\longrightarrow}\pi, \text{ in }L^\infty(\R^+) \quad \text{  and } \quad u_1(\cdot,t)\underset{t\to\infty}{\longrightarrow}-\pi \text{  in } L^\infty_\mathrm{loc}(\R^+);\]
\item unique in the sense that 
the family $\{u_1(\cdot,\cdot+t_y)\}$ contains all entire solutions of (\ref{pb1})--$(\mathcal{BC}_1)$ between $-\pi$ and $\pi$.
\end{enumerate}
Moreover, the homoclinic $u_1$ is the limit of the relative time-periodic solutions  of (\ref{pb1})-(\ref{def:gtheta}) from Theorem \ref{thm1}, parameterized as  $u(x,t;\vartheta,y)$ 
with the convention  $u(0,0;\vartheta,y)=y$, in the sense that 
\begin{enumerate}\setcounter{enumi}{2}
 \item 
For all $y\in(-\pi,\pi)$, there exists $t_y\in\R$ such that 
 \[
 \displaystyle u(\cdot,\cdot\ ; \vartheta,y)\underset{\vartheta\nearrow1}{\longrightarrow}u_1(\cdot,\cdot+t_y) \text{ in  } \displaystyle L^\infty_\mathrm{loc}\left(\R^+\times\R\right).
 \]
 Moreover, the map $y\mapsto t_y$ is a decreasing bijection from $(-\pi,\pi)$ onto $\R$. 
 \item For $y=\pi$,
$ u(\cdot,\cdot\ ;\vartheta,y)\underset{\vartheta\nearrow1}{\longrightarrow}\pi$.
\end{enumerate}
\end{Theorem}
 We emphasize that the convergence of the homoclinic in time is not uniform in $x$, since $u_1(x,t)\to\pi$ for $x\to\infty$ and any fixed $t$. 
 
 Increasing $|\vartheta|$  past the boundary $|\vartheta|=1$, we find two constant stationary solutions as the simple zeros of $g(u;\vartheta)$. 
 Heteroclinic solutions for $|\vartheta|>1$ are found in the next result. 
 
 \begin{Theorem}[Heteroclinic connections]\label{thm3}
  Consider (\ref{pb1}) with generic boundary condition $g$ and with $c>0$. 
  Let $y_1<y_2$ be two consecutive zeros of $g$. Then there exists a unique (up to time translation) entire solution $U_\infty=U_\infty(x,t)$
  of (\ref{pb1}) satisfying $y_1<U_\infty<y_2$. Moreover, if $g>0$ on $(y_1,y_2)$, then 
  \begin{equation*}
   U_\infty(\cdot,t)\underset{t\to-\infty}{\longrightarrow}y_2, \text{ in } L^\infty(\R^+)   \ \textrm{ and }   U_\infty(\cdot,t)\underset{t\to\infty}{\longrightarrow}y_1 \text{ in } L^\infty_\mathrm{loc}(\R^+).
  \end{equation*}
If $g<0$ on $(y_1,y_2)$, the same result holds, interchanging $y_1$ and $y_2$.
 \end{Theorem}

 \paragraph{Purely diffusive case $c=0$.}
 We turn to the somewhat more subtle case of vanishing transport, 
 \begin{equation}
 \begin{cases}\label{pb:ceq0}
  \partial_t u = \partial_{xx}u, \qquad & x>0,\ t>0\\
  \partial_x u = g(u), \qquad & x=0,\ t>0.
 \end{cases}
\end{equation}
 Clearly,  linear functions of the form $u(x)=y+g(y)x$ are stationary solution of (\ref{pb:ceq0}). On the other hand, subtracting linear profiles 
 $u\mapsto u-xy$ simply changes $g$ by adding a constant. We can therefore study dynamics in general by restricting to  solutions bounded in $x$,  
 that is, in  $Y_\mathrm{c}$ with $c=0$, and studying the different cases of $g>0$, $g\geq0$ and sign-changing $g$ as in the previous section. 

In the case $g>0$, we still expect decaying solutions, but diffusive transport is too weak to accommodate relative time-periodic dynamics, but 
merely creates a diffusive drift comparable to the caricature case $g\equiv 1$, as the following result shows. 

 \begin{Proposition}[Diffusive drift]\label{prop:4}
  Assume that $g>0$. Then there is a constant $C>0$ such that, for any bounded initial condition $u_0\in X_+$, the solution $u(x,t)$ of the Cauchy problem (\ref{pb:ceq0})  with $u(\cdot,0)=u_0$ satisfies
  \begin{equation*}
-C\sqrt{t}-\left\Vert u_0\right\Vert_\infty\leq u(0,t)\leq -\frac{\sqrt{t}}{C}+\left\Vert u_0\right\Vert_\infty.
  \end{equation*}
 \end{Proposition}

 In the case where $g$ is not strictly positive, we can still establish the existence of homoclinic and heteroclinic solutions but uniqueness results are weaker.

 \begin{Proposition}[Heteroclinic connections]\label{p:het0}
  Assume that for $y_1<y_2$, we have $g(y_1)=g(y_2)=0$ and $g(u)>0$ for $u\in (y_1,y_2)$. Then there exists a heteroclinic solution $U_\infty(t,x)$ with 
   \begin{equation*}
   U_\infty(\cdot,t)\underset{t\to-\infty}{\longrightarrow}y_2,\text{ in }L^\infty(\R^+) \quad \textrm{ and }   
   \quad U_\infty(\cdot,t)\underset{t\to\infty}{\longrightarrow}y_1 \text{ in }L^\infty_\mathrm{loc}(\R^+).
  \end{equation*}

Moreover, any entire solution 
 between $y_1$ and $y_2$ must be a heteroclinic connection in this sense. 
 \end{Proposition}
Unfortunately, we cannot prove uniqueness of the 
 entire solution. We will briefly discuss some results towards uniqueness in the last section. 
 
 \paragraph{Techniques.}
 Our result rely heavily on comparison principles and generally monotonicity methods. Existence of periodic orbits is established using a homotopy
 and global continuation results, invoking a Krein-Rutman argument as a key ingredient. Heteroclinic and homoclinic orbits are found as suitable limits
 of trajectories, exploiting local compactness. Uniqueness results for heteroclinic orbits are based on invariant manifolds theory. 
 
 \paragraph{Beyond phase-diffusion: boundary reactions and discrete gauge symmetries.}
Our interest here is somewhat specific and we are not aware of results in the literature addressing these specific questions. On the other hand, 
there are numerous results in the literature addressing related problems. An important class of related problems arises for instance in combustion, 
where substances diffuse in a domain but reaction is localized on the boundary, an effect which can be modeled for instance through a nonlinear flux 
as described here; see for instance 
 \cite{Cabre_Consul_Mande,caffarelli}. In the absence of advection, such problems are intrinsically connected with fractional differential equations. 
 
 Our pursuit of non-convergent, periodic dynamics in scalar equations was preceded by 
  \cite{Guidotti_Merino} where periodic solutions in a heat equation with nonlinear, nonlocal boundary conditions are found. 
  
Our specific choice of gauge symmetry on the other hand mimics models for phase oscillators, prominent in particular in neuroscience  
\cite{ermentrout}, where $u\in \R\mod2\pi\Z$ would describe the phase of an individual neuron, and $x$ would parameterize a collection 
of neurons. In this context, the case of $g>0$ is usually referred to as a bursting or oscillatory state, while sign-changing $g$ models 
excitable states; see also \cite{bellay} for the transition between these two scenarios. We emphasize however that the kinetics in this 
scenario are present in the entire domain. 
 
Closer in spirit to our motivation are results in  \cite{nakamura3,nakamura2}, where monotonically increasing solutions in a scalar 
reaction-diffusion system are understood as periodic orbits relative to a gauge symmetry,  similar to our situation. Motivation here 
is contextually related to ours, originating from crystal growth. 

 \paragraph{Outline.} We first recall well-posedness results, comparison principles, and lap number monotonicity results for our 
 context in Section \ref{s:2}. We prove Theorem \ref{thm1} in Section \ref{section:3}. We study the limit $|\vartheta|\to 1$, and,
 more generally, the dynamics between constant solutions in Section \ref{s:4}. Most of the results here apply to the case $c=0$, 
 as well. Section \ref{s:5} presents some 
 specific aspects of the diffusive case $c=0$. We conclude with a discussion.

 \section{Preliminaries}\label{s:2}
 \paragraph{Symmetrization and well-posedness of the Cauchy problems.}
For any $u_0\in Y_\mathrm{c}$, we consider the Cauchy problem (\ref{pb1}),
\begin{equation}
  \begin{cases}\label{CauchyPb}
  \partial_t u = \partial_{xx}u-c\partial_x u, \qquad & x>0,\ t>0\\
  \partial_x u = g(u;\vartheta), \qquad & x=0,\ t>0, \\
  u(x,0)=u_0(x),\qquad & x\geq0.
 \end{cases}
\end{equation}
We claim that there exists a unique classical solution $u=u(x,t;u_0)$ of (\ref{CauchyPb}) such that $u(x,t)\to u_0(x)$ as $t\to0^+$ and $u(\cdot,t)\in Y_\mathrm{c}$, for all $t>0$.
Although the boundary condition in (\ref{CauchyPb}) is less common,  global well-posedness follows from abstract results on semi-linear parabolic equations; 
see for instance \cite{Amann_93,Henry}, or Chapter XIII in \cite{Lieberman}. For a more extensive bibliography on the subject, see \cite{Latushkin_JEE06} and references therein. For the sake of completeness we outline a more self-contained approach, here. We make the dependence on the parameter
$\vartheta$ as in (\ref{def:gtheta}) explicit for later use. Generally, Lipschitz dependence of $g$ on a parameter is sufficient for the results stated here. For any $u\in Y_\mathrm{c}$, let
\begin{equation}\label{symmetrization}
 \tilde{u}(x)=\rme^{-\frac{c}{2}x}u(x)\in X_+.
\end{equation}
Then $u(x,t)$ is a classical solution of (\ref{CauchyPb}) if and only if $\tilde{u}(x,t)$ is a classical solution of 
\begin{equation}
  \begin{cases}\label{Cauchytilde}
  \partial_t \tilde{u} = \partial_{xx}\tilde{u}-\frac{c^2}{4}\tilde{u}, \qquad & x>0,\ t>0\\
  \partial_x \tilde{u} = g(\tilde{u};\vartheta)-\frac{c}{2}\tilde{u}, \qquad & x=0,\ t>0, \\
  \tilde{u}(x,0)=\rme^{-\frac{c}{2}x}u_0(x),\qquad & x\geq0,
 \end{cases}
\end{equation}
and well posedness of (\ref{CauchyPb}) in $Y_\mathrm{c}$ is equivalent to well-posedness of (\ref{Cauchytilde}) in $X_+$. 
The method now relies on standard considerations for the heat equation on the 
half-line (see \cite{Cannon}, chapter 4). Since the linear operator in (\ref{Cauchytilde})
is symmetric, we extend the equation on the full line and consider the boundary 
condition as a localized source term.
Let us consider the embedding $\displaystyle j:X_+\hookrightarrow X=BUC(\R)$ through $j(\tilde{u})(x)=\tilde{u}(|x|)$. For the sake of notation, 
we later identify $j(\tilde{u})$ and $\tilde{u}$ when there is no ambiguity. 
Then, at least formally, $\tilde{u}$ is a solution of (\ref{Cauchytilde}) if and only if 
$j(\tilde{u})$ is a solution of 
$$
\partial_tv = \partial_{xx}v-\frac{c^2}{4}v-2\left\langle\delta_0,g(v)-\frac{c}{2}v\right\rangle\delta_0 
$$
where $\delta_0$ is the Dirac distribution at $x=0.$
All solutions arise as fixed points from a variation-of-constant formula, which leads us to define, for all $T>0$,
\begin{equation}\label{eq:defG}
 G:\left\{\begin{array}{ccl}
           C\left([0,T],X\right)\times X_+\times\R & \to & C\left([0,T],X\right) \\
           (v;\tilde{u}_0,\vartheta) & \mapsto & t\mapsto\Gamma(t)\ast j(\tilde{u}_0)\\
           &&\ \quad\quad -2 \int_0^t\Gamma(\cdot,t-s)\left( g(v(0,s))-\frac{c}{2}v(0,s)\right) \rmd s
          \end{array}
\right.
\end{equation}
where $\displaystyle\Gamma(x,t):=(4\pi t)^{-1/2}\rme^{-\frac{c^2}{4}t-\frac{x^2}{4t}}$ 
is the Green function associated to the first equation of (\ref{Cauchytilde}) on the whole line.
The corresponding semi-group is analytic, hence $G$ is smooth with respect to $\tilde{u}_0$, and the regularity with respect to $\vartheta$ is determined by $g$. 
Any classical solution of (\ref{Cauchytilde})
must be a fixed point of $G$. A direct computation shows that if $\displaystyle T<\pi/ 4\left\Vert\partial_1g(\cdot\ ;\vartheta)\right\Vert_\infty$, then $G$ is a contraction with respect to $v$, 
uniformly in $\tilde{u}_0$, locally uniformly in $\vartheta$, which guarantees the existence and uniqueness of a global solution of 
(\ref{Cauchytilde}) in $X_+$ as desired. As a fixed point of $G$, the solution inherits regularity with respect to $\vartheta$ and $\tilde{u}_0$. 
Standard space-time parabolic regularity shows that the solution is classical.

\paragraph{Sub and super-solutions.}
Throughout, we define a \emph{sub-solution} of (\ref{CauchyPb}) as a function $\underline{u}(x,t)$, 
continuous on $s\leq t\leq T$, and satisfying the following inequalities in the classical sense for all $s<T$,
\begin{equation}
  \begin{cases}\label{subsol}
  \partial_t \underline{u} \leq \partial_{xx}\underline{u}-c\partial_x \underline{u}, \qquad & x>0,\ t\in(s,T)\\
  \partial_x \underline{u} \geq g(\underline{u}), \qquad & x=0,\ t\in(s,T),
 \end{cases}
\end{equation}
We define super-solutions in the same fashion, as functions satisfying (\ref{subsol}) with the reverse inequalities. 
The following classical result is our basic tool; see for instance \cite{smoller} for proofs. 
\begin{Lemma}[Comparison principles]\label{comparison}
 Let $\underline{u}$ be a sub-solution and $\overline{u}$ be a super-solution of (\ref{CauchyPb}), defined for $t\in(s,T)$, such that
 $\underline{u}(x,s)\leq\overline{u}(x,s)$, for all $x\geq0$. Then, either $\underline{u}<\overline{u}$ for all $t\in(s,T)$ or there exists $t_0>s$ such that 
 $\underline{u}\equiv\overline{u}$ for all $t\in[s,t_0)$.
 
 The result remains valid if $\underline{u}$ is a supremum of a finite number of sub-solutions, or if $\overline{u}$ is an infimum of a finite number of super-solutions.
\end{Lemma}
 
 \paragraph{Zero number for parabolic equations.}
 We now turn to  solutions of a linear parabolic equation 
 \begin{equation}\label{eqlin}
 v_t=v_{xx}+c(x,t)v,\qquad x\in I,\ t\in\left( s,T\right),
 \end{equation}
 where $-\infty\leq s<T\leq \infty$, $I=(a,b)$  with $-\infty\leq a < b\leq\infty $ and $c$ is a bounded measurable
 function. We denote by $z_I(v(\cdot,t))$ 
 the number, possibly infinite, of zeros $x\in I$ of the function
 $x\mapsto v(x,t)$. 
 \begin{Lemma}[Zero numbers  \cite{Angenent_Crelle88}]\label{lemzero}
 	Let $v$ be a nontrivial solution of \emph{(\ref{eqlin})} and
 	$I=(a,b)$, with $-\infty\leq a < b\leq \infty$. Assume that the
 	following conditions are satisfied: 
 	\begin{itemize}
 		\item if $b<\infty$, then $v(b,t)\neq0$ for all $t\in\left( s,T\right)$,
 		\item if $a>-\infty$, then $v(a,t)\neq0$ for all $t\in\left( s,T\right)$.
 	\end{itemize}
 	Then the following statements hold true.
 	\begin{enumerate}
 		\item[(i)] For each $t\in\left( s,T\right)$, all zeros of $v(\cdot,t)$ are
 		isolated. In particular, if $I$ is bounded, then
 		$z_I(v(\cdot,t))<\infty$ for all $t\in\left( s,T\right)$. 
 		\item[(ii)] The function $t\mapsto z_I(v(\cdot,t))$ is monotone
 		non-increasing on $(s,T)$ with values in
 		$\N\cup\{0\}\cup\{\infty\}$. 
 		\item[(iii)] If, for some $t_0\in(s,T),$ the function $v(\cdot,t_0)$
 		has a multiple zero in $I$ and $z_I(v(\cdot,t_0))<\infty$, then for
 		any $t_1,t_2\in(s,T)$ with $t_1<t_0<t_2$, one has 
 		\begin{equation}\label{zerodrop}
 		z_I(v(\cdot,t_1))>z_I(v(\cdot,t_0))\ge z_I(v(\cdot,t_2)).
 		\end{equation}
 	\end{enumerate}
 \end{Lemma}
 If (\ref{zerodrop}) holds,  we say that $z_I(v(\cdot,t))$ drops in the
 interval $(t_1,t_2)$. 
 
 \begin{Remark}\label{convzero}{\rm
    If the assumptions of Lemma \ref{lemzero} are
    satisfied and for some $t_0\in(s,T)$, one has
    $z_I(v(\cdot,t_0))<\infty$, and 
    $z_I(v(\cdot,t))$ can drop at most finitely many times in $(t_0,T)$. If $z_I$ is constant on $(t_0,T)$, then $v(\cdot,t)$ has only simple zeros in $I$  
    for all $t\in (t_0,T)$. In particular, if $T=\infty$, there exists
    $t_1<\infty$ such that $t\mapsto z_I(v(\cdot,t))$ is constant on
    $(t_1,\infty)$ and all zeros  
    are simple.	
 }
 \end{Remark}
 
 Using the previous remark and the implicit function 
 theorem, we obtain the following corollary.
 \begin{Corollary}\label{zeroIFT}
 	Assume that the assumptions of Lemma \ref{lemzero} are satisfied and that the function $t\mapsto z_I(v(\cdot,t))$ is constant on $(s,T)$.
 	If for some $(x_0,t_0)\in I\times(s,T)$ one has  $v(x_0,t_0)=0$, then
 	there exists a $C^1$-function $t\mapsto\eta(t)$ defined for
 	$t\in(s,T)$ such that  $\eta(t_0)=x_0$ 
 	and $v(\eta(t),t)=0$ for all $t\in(s,T)$. 
 \end{Corollary}

 \section{Existence of a periodic solution}\label{section:3}
 We prove the existence of the entire solution from  Theorem \ref{thm1}. Recall that $g(\cdot;\vartheta)$  is given through (\ref{def:gtheta}). We shall use continuation in $\vartheta.$
Let $I\subset\left(-1,1\right)$ be the set of all $\vartheta\in(-1,1)$ such that there exist $u=u(x,t;\vartheta)$ and $T=T(\vartheta)>0$ such that $u$ is 
\begin{itemize}
 \item[(P1)] a classical solution of (\ref{pb1});
 \item[(P2)] periodic relative to the gauge: $\displaystyle u\left(\cdot,T(\vartheta)\right) = u\left(\cdot,0\right)-2\pi;$
 \item[(P3)] monotone, $\partial_t u<0$.
\end{itemize}
In order to complete the proof of Theorem \ref{thm1}, we need to prove that $I$ is non-empty and both open and closed in $(-1,1)$.

First, $I$ is nonempty since $0\in I$. Indeed,
\begin{equation}\label{eq:1}
 u(x,t;0)=x-ct
\end{equation}
is a trivial solution of (\ref{pb1})--$(\mathcal{BC}_0)$ satisfying (P1)-(P3), with $T(0)=\frac{2\pi}{c}$. 
Notice also that up to a time shift $u(\cdot,\cdot;0)$ is the only solution of (\ref{pb1})--$(\mathcal{BC}_0)$ 
satisfying (P1)-(P2).
\subsection{The set $I$ is open}
Suppose throughout  that there exists $\vartheta^*\in I$ such that (P1)-(P3) hold true and denote by $u^*(x,t)$ and $T^*$ the corresponding solution and period.

\begin{Proposition}\label{prop:open}
 There exists $\varepsilon>0$ such that $\displaystyle\left(\vartheta^*-\varepsilon,\vartheta^*+\varepsilon\right)\subset I$.
\end{Proposition}
The proof proceeds in several steps.
\paragraph{Relative periodic orbits as fixed points.}
Using the transformation (\ref{symmetrization}), 
the desired periodic property (P2) is  equivalent to 
\begin{equation}\label{periodic:tilde}
 \tilde{u}(x,T) = \tilde{u}(x,0)-2\pi \rme^{-\frac{c}{2}x}.
\end{equation}

We define the map 
\begin{equation}\label{eq:2}
 F : \left\{ \begin{array}{ccl}
              X_+\times(-1,1)\times\R^+ & \longrightarrow & X_+ \\
              \left( \tilde{u}_0,\vartheta,T\right) & \longmapsto & \tilde{u}(\cdot,T)-\tilde{u}_0+2\pi \rme^{-\frac{c}{2}\cdot}
             \end{array}
\right.
\end{equation}
where $\tilde{u}(\cdot,t)$ is the solution of (\ref{Cauchytilde}) with initial condition $\tilde{u}_0$. 
Clearly, zeros of $F$ correspond to relative time-periodic solutions. By assumption, we have that 
\begin{equation}\label{eq:3}
\displaystyle F\left( \tilde{u}^*(\cdot,0;\vartheta^*),\vartheta^*,T^*\right)=0.
\end{equation}

\paragraph{The linearized problem.}
The well-posedness of (\ref{Cauchytilde}) ensures that $F$ is well-defined. Its regularity is given by the regularity of $G$ in the $(\tilde{u}_0,\vartheta)$-variables,
and by the parabolic regularity in $T$. Therefore it is continuously differentiable in all three variables.
Let us denote $\displaystyle D_{\tilde{u}}F^*:=D_{\tilde{u}}F(\tilde{u}^*(\cdot,0,\vartheta^*),\vartheta^*,T^*)$ the partial derivative of $F$ 
with respect to the initial condition at the critical triplet, which is,
\begin{equation}\label{def:DuF}
 D_{\tilde{u}}F^*=\Phi(T^*)-\mathrm{id}
\end{equation}
where $\Phi$ is defined by $\Phi(t)v_0:=v(t)$, where $v(\cdot,t)$ is the solution of the linearized equation
\begin{equation}\label{eq:lin}
 \begin{cases}
   \partial_t v = \partial_{xx}v-\frac{c^2}{4}v, & \qquad t>0,\ x>0 \\
   \partial_x v = \left(\partial_\mathrm{u}g( \tilde{u}^*;\vartheta^*)-\frac{c}{2}\right) v, & \qquad t>0,\ x=0,
 \end{cases}\tag{$\mathcal{L}\mathcal{P}$}
\end{equation}
with initial condition $v(x,0)=v_0(x)$.
Since $\tilde{u}^*$ satisfies (\ref{periodic:tilde}), the time derivative $\partial_t \tilde{u}^*$ is $T^*-$periodic and is a solution of (\ref{eq:lin}).
Hence, together with assumption (P3), 
it yields
\begin{equation}\label{vstar}
 -\partial_t\tilde{u}^*(x,t):=v^*(x,t)>0,\qquad D_{\tilde{u}}F^*[v^*(\cdot,0)]=0.
\end{equation}
In order to complete the proof of Proposition \ref{prop:open}, we show that $D_{\tilde{u}}F^*$ is a Fredholm operator of index 0 with one-dimensional generalized kernel, that is,
\begin{itemize}
 \item $D_{\tilde{u}}F^*$ is invertible up to a compact operator. 
 \item $\displaystyle \ker\, D_{\tilde{u}}F^*=\ker\, \left( D_{\tilde{u}}F^*\right)^2 =  v^*(\cdot,0)\R;$
\end{itemize}
We will then conclude that $D_{\tilde{u},T}F^*$ is onto and solve $F=0$ with the implicit function theorem for nearby parameter values $\vartheta$ after eliminating the kernel.

\begin{Lemma}\label{lemma5}
 There exist linear operators $\mathcal{L}(T^*)$ and $\mathcal{K}(T^*)$ on $X$ such that $\mathcal{L}(T^*)$ is a contraction, $\mathcal{K}(T^*)$ is compact, and 
 $$
 \Phi(T^*)=\mathcal{L}(T^*)+\mathcal{K}(T^*).
 $$ In particular, $D_{\tilde{u}}F^*$ is Fredholm of index 0. 
\end{Lemma}

\begin{Proof}
 Recall the definition of $\Phi(t)v$, $v\in X_+$, 
 as  the solution of (\ref{eq:lin}) with initial condition $v.$ Slightly abusing notation,  we write $v(x,t)$ for the solution with $v(x,0)=v(x)$. 
 Denoting the heat kernel $\displaystyle \Gamma(x,t)=\frac{1}{\sqrt{4\pi t}}\rme^{-\frac{x^2}{4t}}$
 we have the following solution formula for $t\geq0$,
 \begin{align}
  \left[\Phi(t)v\right](x)=v(x,t) = &  \rme^{-\frac{c^2}{4}t}\int_0^\infty\left( \Gamma(x-y,t)+\Gamma(x+y,t)\right) v(y,0)\rmd y \nonumber \\
   & +2\int_0^t \rme^{-\frac{c^2}{4}(t-s)}\Gamma(x,t-s)\left(\frac{c}{2}-\partial_\mathrm{u}g( \tilde{u}^*(0,s);\vartheta^*)\right) v(0,s) \rmd s. \nonumber
 \end{align}
Let us define, for $t>0$ and $v\in X_+$,
\begin{align}
 \left[\mathcal{L}(t)v\right](x) & := \rme^{-\frac{c^2}{4}t}\int_0^\infty\left( \Gamma(x-y,t)+\Gamma(x+y,t)\right) v(y)\rmd y, \label{eq:14}\\
 \left[\mathcal{K}(t)v\right](x) & := 2\int_0^t \rme^{-\frac{c^2}{4}(t-s)}\Gamma(x,t-s)\left( \frac{c}{2}-\partial_\mathrm{u}g( \tilde{u}^*(0,s);\vartheta^*)\right) v(0,s) \rmd s.\label{eq:15}
\end{align}
Clearly, $\mathcal{K}$ is implicitly defined through the solution $v(x,t)$ of (\ref{eq:lin}), but $\mathcal{L}$ is well defined and linear, and so is $\mathcal{K}$, as 
the difference $\Phi-\mathcal{L}$.

\subparagraph{Step 1: $\mathcal{L}(t)$ is a contraction, for all $t>0$.} We quickly estimate  for all $v\in X_+$,
\begin{equation}\label{eq:16}
 \left\Vert \mathcal{L}(t)v\right\Vert_{L^\infty}\leq \rme^{-\frac{c^2}{4}t}\left\Vert v\right\Vert_{L^\infty}.
\end{equation}

\subparagraph{Step 2: $\mathcal{K}$ is uniformly bounded.} This is a direct consequence of the Gronwall Inequality.
For all $v\in X_+,t>0$, we have
\begin{equation}\label{eq:17}
 \mathcal{K}(t)v = 2\int_0^t \rme^{-\frac{c^2}{4}(t-s)}\Gamma(\cdot,t-s)\left[ \left( \mathcal{K}(s)+\mathcal{L}(s)\right) 
 v\right](0)\left( \vartheta\frac{c}{2}-\partial_\mathrm{u}g( \tilde{u}^*(0,s);\vartheta^*)\right) \rmd s.
\end{equation}
Hence, with (\ref{eq:16}) and $\displaystyle \Gamma(x,t)\leq \frac{1}{\sqrt{4\pi t}}$, we have, denoting $C_1=2\left\Vert\partial_\mathrm{u}g\right\Vert_\infty+c$,
\begin{align}
 \left\Vert \mathcal{K}(t)v\right\Vert_{L^\infty} & \leq C_1\int_0^t\frac{\rme^{-\frac{c^2}{4}(t-s)}}{\sqrt{4\pi(t-s)}}\rme^{-\frac{c^2}{4}s}\left\Vert v\right\Vert_{L^\infty}\rmd s +
 C_1\int_0^t\frac{\rme^{-\frac{c^2}{4}(t-s)}}{\sqrt{4\pi(t-s)}}\left\Vert \mathcal{K}(s)v\right\Vert_{L^\infty}\rmd s  \nonumber\\
  & \leq C_1\frac{t}{\sqrt{4\pi}}\rme^{-\frac{c^2}{4}t}\left\Vert v\right\Vert_{L^\infty}+
  C_1\int_0^t\frac{\rme^{-\frac{c^2}{4}(t-s)}}{\sqrt{4\pi(t-s)}}\left\Vert \mathcal{K}(s)v\right\Vert_{L^\infty}\rmd s \nonumber \\
   & \leq \beta(t)\left\Vert v\right\Vert_{L^\infty}+\int_0^t\gamma(s)\left\Vert \mathcal{K}(s)v\right\Vert_{L^\infty}\rmd s, \label{eq:18}
\end{align}
where $t\mapsto\beta(t)$ is positive and bounded and $s\mapsto\gamma(s)$ is positive and of finite integral on $(0,t)$, uniformly in $t$.
From Gronwall's Lemma, there exists a positive and continuous function $C(t)$ such that for all $t\geq0$, for all $v\in X$,
\begin{equation}\label{eq:19}
 \left\Vert\mathcal{K}(s)v\right\Vert_{L^\infty}\leq C(t)\left\Vert v\right\Vert_{L^\infty}.
\end{equation}

\subparagraph{Step 3: $\mathcal{K}(T^*)$ is compact.} By standard parabolic estimates, the restrictions of $\Phi(T^*)$ and 
$\mathcal{L}(T^*)$ to $\displaystyle BUC\left(\left[0,M\right]\right)$ are both compact, for all $M>0$, and so is the restriction of $\mathcal{K}(T^*)$ to 
$\displaystyle BUC\left(\left[0,M\right]\right)$. For all $x\geq M$, we have, combining (\ref{eq:16}-\ref{eq:18}):
\begin{align}
 \left| \left[ \mathcal{K}(T^*)v\right](x)\right| & \leq C_1\int_0^{T^*}\frac{\rme^{-\frac{c^2}{4}(T^*-s)}}{\sqrt{4\pi(T^*-s)}}
 \rme^{-\frac{x^2}{4(T^*-s)}}\left(\left\Vert\mathcal{L}(s)\right\Vert+\left\Vert\mathcal{K}(s)\right\Vert\right)\left\Vert v\right\Vert_{L^\infty}\rmd s \nonumber\\ 
  & \leq \rme^{-\frac{M^2}{4T^*}}C_1\sqrt{\frac{T^*}{4\pi}}\left(1+C(T^*)\right)\left\Vert v\right\Vert_{L^{\infty}}\underset{M\to\infty}{\longrightarrow}0, \label{eq:20}
\end{align}
which establishes the compactness of $\mathcal{K}(T^*)$ on $BUC(\R^+)$. The linearization $D_{\tilde{u}}F^*$ therefore 
is Fredholm of index 0 as a perturbation of the identity by a contraction and a compact map. 
\end{Proof}

\begin{Lemma}\label{lemma3}
 The kernel of $D_{\tilde{u}}F^*$ is one dimensional, 
 $
 \ker\,D_{\tilde{u}}F^*= v^*(\cdot,0)\R.
 $
\end{Lemma}

\begin{Proof}
 Let $w$ be an element in $\displaystyle  \ker\,D_{\tilde{u}}F^*$. With a slight abuse of notation, we denote  by $w(x,t)$ the solution of (\ref{eq:lin})
 with initial condition $w(x,0)=w(x)$.
 We show that $w$ has to be a multiple of $v^*$. The method, standard in parabolic equations, is to prove that, up to a scalar multiple,
 $w$ and $v^*$ can be ordered, and then to apply the parabolic maximum principle.
 The functions $v^*$ and $w$ are both $T^*$-periodic in time, bounded, differentiable, and uniformly continuous for all time. 
 Decomposing into Fourier series, we find, with $\omega={2\pi}/{T^*}$,
 \begin{equation}\label{eq:4}
  v^*(x,t)=\sum_{k=-\infty}^\infty v_k(x)\rme^{\rmi k\omega t},\qquad w(x,t)=\sum_{k=-\infty}^\infty w_k(x)\rme^{\rmi k\omega t}.
 \end{equation}
Since $v,w$ are solutions of (\ref{eq:lin}), the functions $v_k,w_k$ are solutions of the differential equation
\begin{equation}\label{eq:5}
 y''-\left(\frac{c^2}{4}+\rmi k\omega\right) y =0,\qquad x\geq0.
\end{equation}
Solutions of (\ref{eq:5}) are given by $\displaystyle y(x)=A\rme^{\nu_k^+x}+B\rme^{\nu_k^-x}$, where $\displaystyle \nu_k^\pm=\pm\sqrt{\frac{c^2}{4}+\rmi k\omega}$. Notice that 
\begin{equation}\label{eq:6}
 \Re {\nu_k^\pm}=\pm\frac{c}{2}\sqrt{\frac{1+\sqrt{1+16k^2\omega^2/c^4}}{2}}.
\end{equation}
Considering that $\displaystyle \Re {\nu_k^+}=-\Re {\nu_k^-}>0$, because $v$ and $w$ are bounded 
functions, there are complex numbers $\lambda_k,\mu_k$ such that 
\begin{equation}\label{eq:7}
 v_k(x)=\lambda_k\rme^{\nu_k^-x},\qquad w_k(x)=\mu_k\rme^{\nu_k^-x}.
\end{equation}
Moreover, from (\ref{eq:6}) there exists $\delta>0$ such that for 
all $|k|\geq1$, $\Re {\nu_k^-}<-\frac{c}{2}-2\delta$. Thus, 
the asymptotic behavior of $v^*$ and $w$ as $x\to\infty$ is 
\begin{equation}\label{eq:8}
 v^*(x,t)=\lambda_0\rme^{-\frac{c}{2}x}+\rmO\left( \rme^{-\left(\frac{c}{2}+\delta\right) x}\right),\qquad w(x,t)=\mu_0\rme^{-\frac{c}{2}x}+\rmO\left( \rme^{-\left(\frac{c}{2}+\delta\right) x}\right)
\end{equation}
Combining (\ref{eq:8}) with the fact that $\displaystyle \lambda_0=\frac{1}{T^*}\int_0^{T^*}v^*(0,t)dt>0$, 
there exists $k>0$ such that $\displaystyle kv^*(x,t)>\left| w(x,t)\right|$, for all $x\geq0$, for all $t\in[0,T^*)$. Let us define 
\begin{equation*}
 \rho_0:=\inf\left\{ k>0:kv^*(x,t)>|w(x,t)|,x\geq0,t\in[0,T^*)\right\},\textrm{ and } z(x,t)=\rho_0v^*(x,t)-w(x,t). 
\end{equation*}
Then, up to replacing $w$ by $-w$, the function $z$ is a solution of (\ref{eq:lin}) and satisfies:
\begin{equation*}
 z(x,t)\geq0,\ \inf\left\{ z(x,t):x\geq0,t\in[0,T^*)\right\}=0,\textrm{ and }z(x,\cdot) \textrm{ is }T^*-\textrm{periodic.}
\end{equation*}
Using the time periodicity, we are in one of the three following situations.

\subparagraph{Case 1: contact point at the boundary.} There exists $t_0\in[0,T^*)$ 
such that $z(0,t_0)=0$. The function $z$ reaches its minimum at the boundary. 
By the Hopf Lemma, either $z\equiv0$ or $\partial_x z(0,t_0)>0$. But the later contradicts 
the fact that $z$ is a solution of (\ref{eq:lin}), hence $w\equiv\rho_0 v^*$.

\subparagraph{Case 2: finite contact point.} There exist $t_0\in[0,T^*)$ and $x_0>0$ such that $z(x_0,t_0)=0$. Again, by the parabolic maximum principle, 
this implies that $z\equiv0$ and gives the desired result.

\subparagraph{Case 3: infinite contact point.}  $z(x,t)>0$  for all $x,t$, and there exist 
$t_0\in[0,T^*)$ and a sequence $x_n\to\infty$ such that $\displaystyle z(x_n,t_0)\underset{n\to\infty}{\longrightarrow}0$.
In the asymptotics $x\to\infty$, the behavior of $z(x,t)$ is given by (\ref{eq:8}): 
$\displaystyle z(x,t)=\left( \rho_0\lambda_0-\mu_0\right) \rme^{-\frac{c}{2}x}+\rmO\left( \rme^{-\left(\frac{c}{2}+\delta\right) x}\right)$.
By definition of $\rho_0$, it yields $\rho_0\lambda_0-\mu_0=0$. Let $k_0$ be the first $k$ such that $\rho_0\lambda_{k}-\mu_{k}\neq0$ or $\rho_0\lambda_{-k}-\mu_{-k}\neq0$.
Combining (\ref{eq:8}) and (\ref{eq:7}), there exists $\delta_0>0$ such that the asymptotic behavior of $z(x,t)$ as $x\to\infty$ is given by 
\begin{equation*}
  z(x,t)= \rme^{\rmi k_0\omega t}\left(\rho_0\lambda_{k_0}-\mu_{k_0}\right) \rme^{\nu_{k_0}^-x}+\rme^{-\rmi k_0\omega t}\left(\rho_0\lambda_{-k_0}-\mu_{-k_0}\right) 
  \rme^{\nu_{-k_0}^-x}+\rmO\left( \rme^{\left(\Re {\nu_{k_0}^-}-d_0\right) x}\right).
\end{equation*}
But the left part of the above asymptotic behavior is sign changing over $t$, which contradicts $z(x,t)>0$. Hence, this situation is impossible.
As a result, $\rho_0v^*-w\equiv0$, and the proof of the lemma is complete.
\end{Proof}

\begin{Lemma}\label{lemma4}
The eigenvalue 0 is algebraically simple, 
$
  \ker\,\left( D_{\tilde{u}}F^*\right)^2= \ker\,D_{\tilde{u}}F^*.
$
\end{Lemma}

\begin{Proof}
 Let $w\in \ker\,\left( D_{\tilde{u}}F^*\right)^2$. Then $D_{\tilde{u}}F^*[w]$ is an element of $\ker\,D_{\tilde{u}}F^*$, 
 so by the above lemma, there exists $\lambda\in\R$,
$\displaystyle  D_{\tilde{u}}F^*[w]=\Phi(T^*)w-w=\lambda v^*(\cdot,0)$.
We argue by contradiction. Let us assume that $\lambda\neq0$. Then, up to a scaling of $w$, we can assume that
\begin{equation}\label{eq:9}
 \Phi(T^*)w-w= v^*(\cdot,0).
\end{equation}
Slightly abusing notation again, we write $w(x,t)$ for the solution of (\ref{eq:lin}) with  initial condition $w(x,0)=w(x)$. 
Then, by (\ref{eq:9}) it yields $\displaystyle w(\cdot,T^*)-w(\cdot,0)=v^*(\cdot,0)$.
Now the function 
\begin{equation*}
 z(x,t)=w(x,t)-\frac{t}{T^*}v^*(x,t).
\end{equation*}
is a solution to 
\begin{equation}\label{eq:linpourz}
  \begin{cases}
   \partial_t z = \partial_{xx}z-\frac{c^2}{4}z-\frac{1}{T^*}v^*, & \qquad t>0,\ x>0 \\
   \partial_x z = \left(\partial_\mathrm{u}g( \tilde{u}^*;\vartheta^*)-\frac{c}{2}\right) z, & \qquad t>0,\ x=0.
 \end{cases}\tag{$\mathcal{A}\mathcal{P}$}
\end{equation}
Moreover, $z(\cdot,T^*)=z(\cdot,0)$, so $z$ is $T^*$-periodic, as is $v^*$. Notice also that, since $v^*>0$ and $v^*$ solves (\ref{eq:lin}),
then $kv^*(x,t)$ is a super-solution for (\ref{eq:linpourz}) for all $k\in\R$.
The strategy is, again, to derive a contradiction from the maximum principle. Decomposing $z$ onto Fourier series, still with $\omega={2\pi}/{T^*}$,
\begin{equation*}
 z(x,t)=\sum_{k=-\infty}^\infty z_k(x)\rme^{\rmi k\omega t},
\end{equation*}
and, combining (\ref{eq:4}) and (\ref{eq:linpourz}), the functions $z_k$ are solutions of 
\begin{equation*}
 z_k''-\left( \frac{c^2}{4}+\rmi 
 \omega k\right) z_k=\frac{v_k}{T^*}.
\end{equation*}
Let us recall that $v_k$ is a solution of (\ref{eq:5}) and given by (\ref{eq:7}). Using the boundedness of $z$, we find that for some 
complex numbers $b_k$,
\begin{equation}
 z_k(x)=\frac{\lambda_k}{2\nu_k^-T^*}x\rme^{\nu_k^-x}+b_k\rme^{\nu_k^-x}.
\end{equation}
This, along with (\ref{eq:6}) gives us the asymptotic behavior of $z$ as $x\to\infty,$
\begin{equation}\label{eq:11}
 z(x,t)=-\frac{\lambda_0}{cT^*}x\rme^{-\frac{c}{2} x}+b_0\rme^{-\frac{c}{2} x}+\rmO\left( \rme^{-\left(\frac{c}{2}+\delta\right) x}\right).
\end{equation}
Recall that $\lambda_0=\frac{1}{T^*}\int_0^{T^*}v(0,t)dt>0$. The asymptotics (\ref{eq:11}) imply that there exists $M>0$ such that 
\begin{equation}\label{eq:12}
 z(x,t)<0,\textrm{ for all }x>M,\textrm{ for all }t\in[0,T^*).
\end{equation}
 We consider two separate cases.

\subparagraph{Case 1: $z$ is sign-changing.} Let us assume that $z(x_0,t_0)>0$ for some $x_0,t_0$. Because $v^*(x,t)>0$ is time periodic and continuous
there exists $\eta>0$ such that $v^*(x,t)>\eta$ for all $x\leq M$, for all $t\in[0,T^*)$. Hence, with (\ref{eq:12}), there exists $k>0$ such that 
$kv^*(x,t)>z(x,t)$ for all $x\geq0$, for all $t\in[0,T^*)$. Let us again define 
\begin{equation*}
  \rho_0:=\inf\left\{ k>0:kv^*(x,t)>z(x,t),x\geq0,t\in[0,T^*)\right\}.
\end{equation*}
Since $z(x_0,t_0)>0$, we have that $\rho_0>0$, and, due to (\ref{eq:12}), there exists a finite contact point 
$(x_1,t_1)$ such that $\rho_0v^*(x_1,t_1)=z(x_1,t_1)$ and $\rho_0v^*(x,t)\geq z(x,t)$ for all $x,t$. But this is impossible, because $\rho_0v^*$ 
is a super-solution for (\ref{eq:linpourz}).

\subparagraph{Case 2: $z$ is negative.} Let us assume that $z(x,t)<0$, for all $x,t$. Comparing (\ref{eq:8}) with (\ref{eq:12}), for all $k>0$, 
there exists $M=M(k)$ such that $\displaystyle\left| kv^*(x,t)\right| <\left| z(x,t)\right|$, for all $x>M$. As a result, 
$\displaystyle \rho_0:=\inf\left\{ k>0:kv^*(x,t)>z(x,t),x\geq0,t\in[0,T^*)\right
\}$ exists and is negative. Moreover, for the same reason,
the contact point cannot be at infinity, so there is a finite contact point, and we obtain the same contradiction. As a consequence, 
our assumption (\ref{eq:9}) is impossible, and the proof of Lemma \ref{lemma4} is complete.
\end{Proof}

\paragraph{Conclusion of the proof of Proposition \ref{prop:open}.}
The conclusion is classical for periodic autonomous systems, showing that periodic orbits persist if the trivial Floquet multiplier is algebraically simple. 
Specifically,  Lemma \ref{lemma5} and \ref{lemma3} imply that the differential $D_{\tilde{u}}F^*$ defined by (\ref{def:DuF})  
has closed range and one-dimensional co-kernel. The derivative of $F^*$ with respect to $T^*$ is $\partial_t\tilde{u}_0^*$ and, by Lemma 
\ref{lemma4}, does not belong to the range of $D_{\tilde{u}}F^*$, such that  $D_{\tilde{u},T}F^*$ is onto. We can now apply the implicit function theorem, 
solving for $\tilde{u}$ in a complement of the kernel. Note that the implicit function theorem also implies local uniqueness. This proves the 
lemma once we establish that $\partial_t u<0$ for the solution. 

In order to show that the condition  $\partial_t \tilde{u}<0$ is  open, notice that $\partial_t \tilde{u}(t,x)<0$ for $x$ sufficiently large, 
due to (\ref{eq:8}), such that uniform continuity on compact sets suffices to establish the claim. 
\qed

\subsection{The set $I$ is closed in $(-1,1)$}
Let $\displaystyle \left(\vartheta_n\right)_n\subset I$ be a sequence satisfying $\displaystyle \vartheta_n\underset{n\to\infty}{\longrightarrow}\vartheta_\infty\in I$. 
Let $T_n,u_n(x,t)$ be the corresponding periods and entire solutions satisfying (P1)-(P3). Without loss of generality, we assume 
\begin{equation}\label{eq:22}
0\leq \vartheta_n\underset{n\to\infty}{\nearrow}\vartheta_\infty,\ \vartheta_\infty>0,\ \textrm{ and }u_n(0,0)=0.
\end{equation}
We need to prove that $\displaystyle \vartheta_\infty\in I$, that is,  there is a period $0<T_\infty<\infty$ and a relative time-periodic entire solution $u_\infty$ of (\ref{pb1})--$(\mathcal{BC}_{\vartheta_\infty})$
satisfying (P1)-(P3). We first derive a priori estimates for $\partial_x u_n$, which will allow us to construct sub and super-solutions to control the behavior of $T_n$.
The conclusion then follows from parabolic regularity.

\begin{Lemma}\label{lemma6}
 For all $n$, we have,
$1-\vartheta_\infty\leq \partial_x u_n(x,t)\leq 1+\vartheta_\infty$ for all $x\geq0,t\in\R$.  
\end{Lemma}

\begin{Proof}
Fix $n$, and let us define $v(x,t)=\partial_x u_n(x,t)$. Then $v$ is a $T_n$-periodic function in time, solution of
\begin{equation}\label{eq:23}
  \begin{cases}
  \partial_t v = \partial_{xx}v-c\partial_x v, \qquad & x>0,\ t\in\R\\
  v(0,t) = f(t), \qquad & x=0,\ t\in\R
 \end{cases}
\end{equation}
with $\displaystyle f(t)=1+\vartheta h(u_n(0,t))$. As a result, $v(x,\cdot)$ is periodic and at the boundary we have
\begin{equation}\label{eq:24}
 1-\vartheta_\infty\leq v(0,t)\leq 1+\vartheta_\infty,\ t\in\R.
\end{equation}
Moreover, by parabolic estimates, considering that $u_n\in Y_\mathrm{c}$, for all $\delta>0$, there exists $C_1>0$ 
such that 
\begin{equation}\label{eq:25}
 \left| v(x,t)\right|\leq C_1\rme^{\left(\frac{x}{2}+\delta\right) x},\qquad x\geq0,t\in\R.
\end{equation}

\textbf{Step 1.} We claim that for all $\delta>0$, there exists $C_2=C_2(\delta)$ such that 
\begin{equation}\label{eq:26}
  \left| v(x,t)\right|\leq C_2\rme^{\delta x},\qquad x\geq0,t\in\R.
\end{equation}
To prove this claim, decompose $v(x,t)$ onto Fourier series:
\begin{equation*}
 v(x,t)=\sum_{k=-\infty}^\infty v_k(x)\rme^{\rmi k\omega t},\qquad \omega=\frac{2\pi}{T_n}.
\end{equation*}
Then, from (\ref{eq:23}), the functions $v_k$ are solutions of the differential equation
\begin{equation}\label{eq:27}
 v_k''-cv_k'-\rmi k\omega v_k=0,\ x\geq0.
\end{equation}
Solutions of (\ref{eq:27}) are 
\begin{equation*}
 v_k(x) = A_k\rme^{\eta_k^+x}+B_k\rme^{\eta_k^-x},\ \eta_k^\pm=\frac{c\pm\sqrt{c^2+4\rmi k\omega}}{2}.
\end{equation*}
This gives $\displaystyle \Re {\eta_k^\pm}=\frac{c}{2}\pm\frac{1}{2}\Re {\sqrt{c^2+4\rmi k \omega}}$, and, in particular,
$\displaystyle \Re {\eta_k^+}>c$ and $\Re {\eta_k^-}<0$, for all $k$. Taking into account (\ref{eq:25}), the first inequality implies 
$A_k=0$, for all $k$, and, together with the second, it implies the desired estimate (\ref{eq:26}).

\textbf{Step 2. }We conclude the proof of Lemma \ref{lemma6} by a sliding argument. 
Let us define, for parameters $k,\eta>0$:
\begin{equation}
 \overline{v}(x;k,\eta)=1+\vartheta_\infty+k\rme^{\eta x}, \ \underline{v}(x;k,\eta)=1-\vartheta_\infty-k\rme^{\eta x}.
\end{equation}
Then, for all $0<\eta<c$, $\overline{v}$ and $\underline{v}$ are super and sub-solution 
of (\ref{eq:23}) respectively. We prove the upper bound of Lemma \ref{lemma6}, the lower being
similar. Fix $\eta>0$. Then, by (\ref{eq:26}), $\overline{v}(\cdot;k,\eta)>v$ for $k$ large enough. 
Let $\displaystyle k_0:=\inf\left\{k\geq0:\overline{v}(x;k,\eta)>v(x,t),x\geq0,t\in\R\right\}$. Then necessarily $k_0=0$.
Otherwise the function $w(x,t):=\overline{v}(x;k_0,\eta)-v(x,t)$ is a time-periodic super-solution of (\ref{eq:23}) satisfying 
$w(0,t)>0$ (because of (\ref{eq:24})), $w(x,t)\geq0$ and $w(x_0,t_0)=0$ for some $x_0,t_0$, which is impossible.
As a result, for all $k>0$, for all $\eta>0$, it reads $\overline{v}(x;k,\eta)>v(x,t)$, for all $x\geq0,t\in\R$.
Passing to the limit $k,\eta\to0$, we obtain the upper bound of Lemma \ref{lemma6}.
 \end{Proof}

\begin{Lemma}\label{lemma7}
 There exist $\lambda,\mu>0$ such that for all $n$,
 \begin{equation*}
  \frac{2\pi}{\mu}\leq T_n\leq\frac{2\pi}{\lambda}.
 \end{equation*}
\end{Lemma}

\begin{Proof}
 The proof is based on the construction of sub and super solutions, that are time decreasing with constant speed.
 Let $\displaystyle K_1=\frac{6\vartheta_\infty}{c(1-\vartheta\infty)}$. Let $f$ be a smooth function with $f(0)=0$ such that,
 \begin{equation*}
  \begin{cases}
 f'(0)=1-\vartheta_\infty\leq f'(x), \textrm{ for all }x\geq0. \\
 f'(x)=1+\vartheta_\infty,\textrm{ for all }x>K_1. \\
 \left\Vert f''\right\Vert_\infty\leq\frac{3\vartheta_\infty}{K_1}.
  \end{cases}
 \end{equation*}
For any $M\in\R$, define
\begin{equation}
\lambda :=\frac{c}{2}\left(1-\vartheta_\infty\right),\qquad \overline{u}(x,t):=-\lambda t+f(x)+M.
\end{equation}
Then $\overline{u}$ is a super-solution for (\ref{pb1}) with ($\mathcal{BC}_{\vartheta_n}$) for all $n$, for all $M$. From Lemma \ref{lemma6} 
and our choice of $f$, for some $M=M_1$ large enough, $\displaystyle\overline{u}(x,0)>u_n(x,0)$, for all $x\geq0$, for all $n$.
Moreover, $\displaystyle \partial_x\overline{u}(0,t)=1-\vartheta_\infty\leq\partial_xu_n(0,t)$, for all $t\geq0$. As a consequence, 
$\displaystyle \overline{u}(x,t)\geq u_n(x,t)$, for all $x,t\geq0$.

Similarly, with $\displaystyle K_2= \frac{6\vartheta_\infty}{c(1+\vartheta_\infty)}$, let $g$ be a smooth function with $g(0)=0$ such that:
 \begin{equation*}
  \begin{cases}
 g'(0)=1+\vartheta_\infty\geq g'(x), \textrm{ for all }x\geq0. \\
 g'(x)=1-\vartheta_\infty,\textrm{ for all }x>K_2. \\
 \left\Vert g''\right\Vert_\infty\leq\frac{3\vartheta_\infty}{K_2}.
  \end{cases}
 \end{equation*}
Then, define for any $M\in\R$,
\begin{equation*}
\mu :=\frac{c}{2}\left(1+\vartheta_\infty\right),\qquad \underline{u}(x,t):=-\mu t+g(x)-M.
\end{equation*}
Then $\underline{u}$ is a sub-solution for (\ref{pb1}) with ($\mathcal{BC}_{\vartheta_n}$) for all $n$, for all $M$. 
From our choice of $g$, for some $M=M_2$ large enough, $\displaystyle\underline{u}(x,0)<u_n(x,0)$, for all $x\geq0$, for all $n$.
Moreover, $\displaystyle \partial_x\underline{u}(0,t)=1+\vartheta_\infty\geq\partial_xu_n(0,t)$, for all $t\geq0$. As a result, we have 
the following inequalities,
\begin{equation}\label{eq:30}
 -\mu t+g(x)-M_2<u_n(x,t)<-\lambda t+f(x)+M_1,\ \textrm{ for all }t\geq0,n\geq0,x\geq0.
\end{equation}
Considering that the functions $u_n$ are relative time-periodic in the sense of (P2), letting $t\to\infty$ in (\ref{eq:30}),
we obtain the desired bounds.
\end{Proof}
\begin{Proof}[ of Theorem \ref{thm1}.] By Lemma \ref{lemma7}, there exists $\displaystyle 0<T_\infty<\infty$
such that, up to a subsequence, $\displaystyle T_n\to T_\infty$. By parabolic regularity, there exists $u_\infty(x,t)$
such that, up to a subsequence, $u_n$ and its derivatives converge to $u_\infty$ and its derivatives, locally uniformly in 
$\displaystyle (x,t)\in\R^+\times[0,T_\infty]$. Then $u_\infty$ solves (\ref{pb1}) with ($\mathcal{BC}_{\vartheta_\infty}$) and satisfies 
$u_\infty(\cdot,T_\infty)=u_\infty(\cdot,0)-2\pi$. Moreover, $\partial_t u_\infty\leq0$ as a limit of positive functions,
so necessarily $\partial_t u_\infty<0$ due to the maximum principle. Hence $\vartheta_\infty$ belongs to $I$, so $I=(-1,1)$.
\end{Proof}

\begin{Remark}
 Considering that (\ref{eq:1}) is the unique solution of (\ref{pb1}) with ($\mathcal{BC}_{0}$), the implicit functions theorem 
 also provides the uniqueness of the solution under the condition $\partial_t u<0$. We do not know if it remains unique without this condition. 
 However, no other entire solutions can be continuously connected to our monotone branch.
\end{Remark}

 \section{Dynamics in the presence of stationary solutions}\label{s:4}
 We now turn to the case when  $g(y)=0$ for some $y\in\R$. 
 We start establishing convergence of the solutions from  Theorem \ref{thm1} in the limit $\vartheta\nearrow 1$, the basis of the proof of 
Theorem \ref{prop2}. We then discuss asymptotics of solutions of (\ref{pb1})  as 
$t\to\pm\infty$ in Section \ref{sub:longtime}. Most of the results there are valid in both cases $c>0$ and $c=0$. We restrict to $c>0$ only 
in the last paragraph, concerned with the uniqueness, and prove  Theorem \ref{thm3}. 
Theorem \ref{prop2} is then  a direct consequence.
 
 \subsection{The case $\vartheta=1$ as a limit}\label{sub:theta1}
 We consider boundary conditions $(\mathcal{BC}_\vartheta)$ with the additional assumption $h(y)>-1$ for $y\in(-\pi,\pi)$.
 For all $y\in(-\pi,\pi]$ and for all $\vartheta\in(-1,1)$, let $u(\cdot,\cdot\ ;\vartheta,y)$ be the entire solution of (\ref{pb1}) given by Theorem \ref{thm1} with 
 $u(0,0;\vartheta,y)=y$.
 
 \begin{Lemma}\label{lemma10}
  For all $M>0,\tau>0$ the family $\{ u(\cdot,\cdot\ ;\vartheta,y):\vartheta\in(-1,1)\}$ in uniformly bounded in $\displaystyle L^\infty\left([0,M]\times[-\tau,\tau]\right)$.
 \end{Lemma}

 \begin{Proof}
  This lemma is a direct consequence of Lemmas \ref{lemma6} and \ref{lemma7}.
  In the proof of the latter, by construction we obtain that $T(\vartheta)\geq\frac{2\pi}{c}$ for all $\vartheta\in(-1,1).$ Therefore by $T-$periodicity it yields that for 
  all $t\in[-\tau,\tau],$ $|u(0,t;\vartheta,y)|\leq\tau c+|y|.$ Lemma \ref{lemma6}
ensures that $\displaystyle 0<\partial_x u(\cdot,\cdot\ ;\vartheta,y)<2,$ which completes the proof of boundedness, locally in space and time  around $(x,t)=
(0,0).$
 \end{Proof}

 \paragraph{Convergence of $u(\cdot,\cdot\ ; \vartheta,\pi)$.}

 Let $(\vartheta_n)\subset(-1,1)$ be any sequence with $\vartheta_n\nearrow1$. Then, by Lemma \ref{lemma10} and parabolic regularity, up to a subsequence, 
 the sequence $u(\cdot,\cdot\ ; \vartheta_n,\pi)$ converges to some function $U(x,t)$ in $L^\infty_\mathrm{loc}(\R^+\times \R)$, solution of (\ref{pb1})--$(\mathcal{BC}_1)$.
 Moreover, due to Theorem \ref{thm1} and by definition of $u(\cdot,\cdot\ ; \vartheta_n,\pi)$, it satisfies
 \begin{equation*}
  \partial_tU\leq0,\ \partial_xU\geq0,\ U(0,0)=\pi, \textrm{ hence }U(x,0)\geq\pi.
 \end{equation*}
Since the constant function $\pi$ is a stationary solution of (\ref{pb1})--$(\mathcal{BC}_1)$, the comparison principle implies that $U\equiv\pi$.
Therefore, any sequence $\vartheta_n\nearrow1$ contains a subsequence along which $u(\cdot,\cdot\ ; \vartheta_n,\pi)$ converges to $\pi$ in $L^\infty_\mathrm{loc}(\R^+\times\R)$,
and from Lemma \ref{lemma10} the family $\{u(\cdot,\cdot\ ; \vartheta,\pi):\vartheta\in(-1,1)\}$ is precompact with respect to this topology. Therefore,  convergence 
holds for $\vartheta\nearrow1$ which proves (iv) in  Theorem \ref{prop2}.
\qed

\paragraph{Convergence to some heteroclinic connection between $-\pi$ and $\pi$.}
We apply the same compactness method. Fix any $y_0\in(-\pi,\pi)$ and any sequence $\vartheta_n\nearrow1$. Then, up to a subsequence, the sequence $u(\cdot,\cdot\ ; \vartheta_n,y_0)$
converges to some function $U_0(x,t)$, solution of (\ref{pb1})--$(\mathcal{BC}_1)$ in $L^\infty_\mathrm{loc}(\R^+\times\R)$. 
Taking into account that $u(\cdot,\cdot\ ; \vartheta_n,y_0)<u(\cdot,\cdot\ ; \vartheta_n,\pi)$
and the above convergence result, it satisfies
\begin{equation}\label{eq:40}
 \partial_xU_0\geq0,\ \partial_tU_0\leq0,\ U_0(0,0)=y_0,\ U_0(x,0)\leq\pi.
\end{equation}

Since $\partial_tU_0$ solves a linear equation, either $U_0$ is a stationary solution or $\partial_tU_0<0$. But (\ref{pb1})--$(\mathcal{BC}_1)$ does not admit  any bounded stationary
solutions besides the constant functions $\pi+2k\pi$, $k\in\Z$. Therefore, $\partial_tU_0<0$, for all $t\in\R$. Similarly, 
(\ref{eq:40}) and the comparison principle imply that $\displaystyle -\pi<U_0<\pi$, for all $x,t$. Due to its monotonicity, the function $U_0$ 
converges to some stationary solutions as $t\to\pm\infty$, and the only possibility is 
\begin{equation}
 U_0(\cdot,t)\underset{t\to-\infty}{\longrightarrow}\pi,\qquad U_0(\cdot,t)\underset{t\to\infty}{\longrightarrow}-\pi.
\end{equation}
This proves (iii) in Theorem \ref{prop2}.

\subsection{Dynamics between two stationary solutions}\label{sub:longtime}
This subsection is devoted to the proof of Theorem \ref{thm3} and (i)--(ii) of Theorem \ref{prop2}.  Consider (\ref{pb1}) with any boundary term $g$, and let $y_1>y_2$ be two consecutive zeros of $g$.
Without loss of generality, we assume 
\begin{equation}\label{eq:41}
 y_1>y_2=0,\qquad g(y)<0\ \textrm{ for }y\in(0,y_1),\qquad g(0)=g(y_1)=0.
\end{equation}
Therefore, the constant functions $U_-\equiv0$ and $U_+\equiv y_1$ are stationary solutions of (\ref{pb1}),
 and there is no other stationary solution with range in $[0,y_1]$.
 
 \paragraph{Existence of the heteroclinic connection.}
 
 We construct the entire solution $U_\infty$ in a general setting, using sub and super-solutions.
 The following key ingredient gives a long time control of the super-solution.
 \begin{Lemma}\label{lemma:sup}
  Fix $c\geq 0$. For all large enough $k$, there exist $T_k>0$ and a function 
  $\displaystyle \overline{u}_k\in C\left([0,T_k],X_+\right)$ such that
  \begin{enumerate}
   \item $\overline{u}_k(\cdot,0)=1/k^2$;
   \item the function $\overline{u}_k$ is a super-solution of (\ref{pb1}) for $t\in(0,T_k]$;
   \item $\overline{u}_k(0,T_k)=1/k$;
   \item $T_k\to\infty$ as $k\to\infty$.
  \end{enumerate}
 \end{Lemma}

 \begin{Proof}
  Fix $\varepsilon>0$. There exists $\delta>0$ such that 
  \begin{equation}\label{eq:44-1}
   g(y)\geq(g'(0)-\varepsilon)y,\ \textrm{ for all }y\in[0,\delta].
  \end{equation}
Let us define 
\begin{equation}\label{eq:46-2}
 \lambda:=-\frac{c^2}{4}+4\left(\frac{c}{2}-g'(0)+\varepsilon\right)^2>0,\qquad \gamma:=\sqrt{\frac{c^2}{4}+\lambda}.
\end{equation}
Fix $k>1/\delta$ and define $\overline{u}_k$ as the solution of
\begin{equation*}
  \begin{cases}
  \partial_t \overline{u}_k = \partial_{xx}\overline{u}_k-c\partial_x \overline{u}_k, \qquad & x>0,\ t>0\\
  \partial_x \overline{u_k} = -\frac{\gamma}{k^2}\rme^{\lambda t}, \qquad & x=0,\ t>0, \\
  \overline{u}_k(x,0)=1/k^2,\qquad & x\geq0.
 \end{cases}
\end{equation*}
Then, considering (\ref{eq:44-1}), $\overline{u}_k$ is a super-solution as long as 
$\overline{u}_k(0,t)\leq\delta$. Therefore, since $k>1/\delta$, we define
\begin{equation*}
 T_k:=\inf\left\{t>0: \overline{u}_k(0,t)=\frac{1}{k}\right\}.
\end{equation*}
The proof of Lemma \ref{lemma:sup} is complete if we can prove that $T_k$ is well defined 
and diverges as $k\to\infty$. The linear form of (\ref{eq:defG}) gives
\begin{equation*}
 \overline{u}_k(0,t)=\frac{\rme^{-\frac{c^2}{4}t}}{k^2\sqrt{4\pi t}}\int_\R \rme^{-\frac{x^2}{4t}-
 \frac{c}{2}|x|}\rmd x+2\frac{\gamma}{k^2}\int_0^t \frac{\rme^{-\frac{c^2}{4}(t-s)}}{\sqrt{4\pi(t-s)}}\rme^{\lambda s}\rmd s.
\end{equation*}
A direct computation leads to 
\begin{equation}
 \overline{u}_k(0,t)=\frac{1}{k^2}\textrm{erfc}\left(\frac{c\sqrt{t}}{2}\right)+\frac{\rme^{\lambda t}}{k^2}\textrm{erf}\left(\gamma^2\sqrt{t}\right),
\end{equation}
where $\displaystyle \textrm{erf}(x)=\frac{2}{\sqrt{\pi}}\int_0^x\rme^{-u^2}\rmd u$ and $\displaystyle \textrm{erfc}(x)=1-\textrm{erf}(x)$. Therefore, $T_k$ is well defined, and satisfies
\begin{equation}
 T_k>\frac{1}{\lambda}\log(k-1)\underset{k\to\infty}{\longrightarrow}\infty.
\end{equation}

 \end{Proof}

 We now have all the ingredients to construct $U_\infty$. For all integer $n$ large enough,
 let $\underline{u}_n(x,t)$ be the solution of the Cauchy problem (\ref{CauchyPb}) starting 
 from the initial condition 
 \begin{equation}\label{eq:46-3}
  \underline{u}_n(x,0)=\max\left(\frac{1}{n}+g\left(\frac{1}{n}\right) x,0\right).
 \end{equation}
Notice that (\ref{eq:46-3}) defines a sub-solution as the supremum of two sub-solutions, 
therefore the  solution satisfies $\partial_t\underline{u}_n(\cdot,t)>0$ and $\partial_x\underline{u}_n(\cdot,t)<0$ for all $t>0$. As a consequence,  $\underline{u}_n(\cdot,t)$
converges to some stationary solution of (\ref{pb1}) as $t\to\infty$. The only possibility is 
$U_+\equiv y_1$, therefore 
\begin{equation}\label{eq:46-4}
 \underline{u}_n(\cdot,t)\underset{t\to\infty}{\longrightarrow}y_1\ \textrm{ in }L^\infty_\mathrm{loc}(\R^+).
\end{equation}
Moreover, for all $n,k$ with $n>k^2$ we have $\underline{u}_n(x,0)<\overline{u}_k(x,0)$, and therefore
\begin{equation}\label{eq:46-5}
 \underline{u}_n(0,t)<\overline{u}_k(0,t)<\frac{1}{k}, \ \textrm{ for all }t\in(0,T_k).
\end{equation}
From (\ref{eq:46-4}), we deduce that there exists $\tau_n>0$ such that $\displaystyle \underline{u}_n(0,\tau_n)=y_1/2$. From (\ref{eq:46-5}) and Lemma \ref{lemma:sup}, we obtain that $t_n\to\infty$ 
as $n\to\infty$. Let us define the sequence
\begin{equation*}
 u_n(x,t):=\underline{u}_n(x,t+\tau_n),\ \textrm{ for }x\geq0,t>-\tau_n.
\end{equation*}
By parabolic regularity, there exists an extraction $\varphi:\N\upuparrows\N$ such that 
the sequence $(u_{\varphi(n)})$ converges to some function $U_\infty$ together with its derivatives 
in $L^\infty_\mathrm{loc}(\R^+)$, and $U_\infty$ is an entire solution of (\ref{pb1}) for all time 
$t\in\R$. It inherits the following properties from the sequence $u_n$, 
\begin{equation}\label{eq:46-6}
 0\leq U_\infty\leq y_1,\ \partial_tU_\infty\geq0,\ \partial_xU_\infty\leq0,\ U_\infty(0,0)=y_1/2,\ U_\infty(\cdot,t)\underset{t\to\infty}{\longrightarrow}0\textrm{ in }L^\infty(\R^+).
\end{equation}
Since $U_\infty$ cannot be a stationary solution, the inequalities in (\ref{eq:46-6}) are strict 
inequalities, and $U_\infty(\cdot,t)$ is convergent to some stationary solution as $t\to\infty$.
The only possibility is $U_+\equiv y_1$.

We have thus established the following result. 

\begin{Proposition}\label{prop:17}
 Consider (\ref{pb1}) with any $c\geq0$. Assuming (\ref{eq:41}) for the nonlinear term $g$,
 there exists an entire solution $U_\infty(x,t)$ of (\ref{pb1}) such that
 \begin{equation*}
 U_\infty(\cdot,t)\underset{t\to-\infty}{\longrightarrow}0\ \textrm{ in }L^\infty(\R^+),\qquad 
  U_\infty(\cdot,t)\underset{t\to\infty}{\longrightarrow}y_1\ \textrm{ in }L^\infty_\mathrm{loc}(\R^+).
 \end{equation*}
\end{Proposition}
This results includes the existence part of Theorem \ref{thm3}. Along with the uniqueness part of Theorem \ref{thm3}, which will be proved in the next section, 
it also concludes the proof of Theorem \ref{prop2}.

\paragraph{Long-time behavior of solutions.}
We consider  solutions of (\ref{pb1}) with $g$ satisfying (\ref{eq:41}).
Let us recall some standard facts about limit sets of a solution.

For any forward solution $u(x,t)$ of (\ref{pb1}), uniformly bounded in $t\geq0$, we define 
its $\omega$-limit set
\begin{equation}\label{def:omegaset}
	\omega(u):=\left\{ \varphi:u(\cdot,t_n)\to\varphi\textrm{ for some sequence }t_n\to\infty\right\}.
\end{equation}
 Similarly, for a bounded ancient solution $U(x,t)$, defined for all $t\leq 0$,  of (\ref{pb1}), we defined its $\alpha$-limit 
 set
 \begin{equation}\label{def:alphaset}
	\alpha(U):=\left\{ \varphi:U(\cdot,t_n)\to\varphi\textrm{ for some sequence }t_n\to-\infty\right\}.
\end{equation}
 In both cases, the convergence is understood in $L^\infty_\mathrm{loc}(\R^+)$.
 As a consequence of parabolic estimates, the $\omega$ and
$\alpha$-limit set of a bounded solution are non-empty, compact and connected for the 
considered topology. Moreover, they are invariant with respect to the evolution problem (\ref{pb1}): 
if $\varphi\in\omega(u),$ then there exists a function $\tilde{U}(x,t)$ entire solution of (\ref{pb1}) such that $\tilde{U}(\cdot,0)=\varphi$ and 
$\tilde{U}(\cdot,t)\in\omega(u)$ for all $t\in\R$. Let us briefly recall how such an entire solution is found. 
If $U(\cdot,t_n)\to\varphi$, consider the sequence $U_n(x,t)=u(x,t+t_n)$. 
Passing to a subsequence if necessary, it converges to $\tilde{U}$ as $n\to\infty$ in 
$C^{2,1}_\mathrm{loc}(\R^+\times\R)$, which therefore also solves (\ref{pb1}). The same holds true 
for an element in $\alpha(U).$

A function $U$ is convergent as $t\to\infty$ (resp. $t\to-\infty$) if its 
$\omega$-limit set (resp. its $\alpha$-limit set) is reduced to a single element.
In this situation, the convergence holds in $C^2_\mathrm{loc}(\R^+)$. Notice also that due to the above argument, it converges necessarily  to a stationary solution. The main result of 
this paragraph is the following proposition, describing the long time behavior of the 
dynamics between two stationary solutions.

\begin{Proposition}\label{prop:18}
 Fix $c\geq 0$. Consider our advection-diffusion equation (\ref{pb1}) with boundary nonlinearity having two zeros $y_1>y_2=0$, $g<0$ on $(y_2,y_1)$, as in (\ref{eq:41}).
\begin{enumerate}
 \item Let $U(x,t)$ be an entire solution with $0\leq U\leq y_1$. Then either 
 $\alpha(U)=\{0\}$ and $\omega(U)=\{y_1\}$, or $U$ is constant.
 \item Let $u(x,t)$ be a solution with initial condition $0\leq u_0\leq y_1$. Then
 either $\omega(u)=\{y_1\}$ or $\omega(u)=\{0\}$.
\end{enumerate}
\end{Proposition}
\begin{Proof}
We first give the proof of the first part of Proposition \ref{prop:18} in the case $c>0$, and 
explain thereafter how to adapt.  Let $U(x,t)$ 
be an entire solution of (\ref{pb1}) such that $0\leq U\leq y_1$.
Recall that the only admissible bounded stationary solutions 
 are the constant functions $0$ and $y_1$. Let us first deal with the $\alpha$-limit set. One of the following mutually exclusive cases
 must hold:
 \begin{itemize}
  \item[(C1)] $\alpha(U)=\{0\}$.
  \item[(C2)] $y_1\in\alpha(U)$.
  \item[(C3)] There exist at least two distinct functions $\varphi_1,\varphi_2\in\alpha(U)$, none of them equal to $y_1$.
 \end{itemize}
We need to prove that the case (C2) implies that $U$ is constant equal to $y_1$, and that 
the case (C3) is impossible.

\subparagraph{Case (C2): $y_1\in\alpha(U)$.}
In this case, by definition, there exists some sequence $t_n\to-\infty$ such that 
\begin{equation}\label{eq:43}
	U(\cdot,t_n)\underset{n\to\infty}{\longrightarrow}y_1\ \textrm{ in }L^\infty_\mathrm{loc}(\R^+).
\end{equation}
For some given $y_0\in(0,y_1)$, let us define
\begin{equation*}
 u_0(x) :=\max\left( y_0+g(y_0)x,0\right).
\end{equation*}
Notice that $u_0$ is the supremum of two sub-solutions, it is therefore a sub-solution 
for the evolution problem. Let $\underline{u}(x,t)$ be the solution of the Cauchy problem (\ref{CauchyPb}) with initial condition $u_0$. 
Since $u_0$ is a sub-solution, $\partial_t\underline{u}>0$ for $t\sim0^+$,
therefore $\partial_t\underline{u}>0$ for all $t>0$, and $\underline{u}(\cdot,t)$ is convergent as $t\to\infty$
in $L^\infty_\mathrm{loc}(\R^+)$ to some stationary solution. The only possibility is $y_1$. Therefore, for all $\varepsilon>0$, 
there exists $T_\varepsilon>0$ such that
\begin{equation}\label{eq:44}
	y_1-\varepsilon<\underline{u}(x,t)<y_1,\ \textrm{ for all }x>\frac{1}{\varepsilon},\ t>T_\varepsilon.
\end{equation} 
By (\ref{eq:43}) there is $N_0$ such that for all $n\geq N_0$, $U(x,t_n)>u_0(x)$. Therefore, the comparison principle
gives
\begin{equation}\label{eq:45}
 U(x,t_n+t)>\underline{u}(x,t) ,\ \textrm{ for all }x>0,\ t>0,\ n\geq N_0.
\end{equation}
Since $t_n\to-\infty$, from (\ref{eq:44}-\ref{eq:45}) necessarily $U\equiv y_1$.

\subparagraph{Case (C3): } there exist two distinct functions $\varphi_1,\varphi_2\in\alpha(U)$, with $\varphi_{1,2}< y_1$.
We first claim that, up to choosing different functions in $\alpha(U)$, we can assume that 
\begin{equation}\label{eq:48}
0\leq \varphi_1(0)<\varphi_2(0)<y_1.
\end{equation}
Indeed, by assumption (C3), one element in $\alpha(U)$ is not equal to $0$,  say $\varphi_1$.
Let $\tilde{U}(x,t)$ be an entire solution of $(\mathcal{P}_1)$ such that $\tilde{U}(\cdot,0)=\varphi_1$ and $\tilde{U}(\cdot,t)\in\alpha(U)$.
Assume by contradiction that it is impossible to satisfy (\ref{eq:48}). Then 
this would imply that $\tilde{U}(0,t)=\varphi_1(0)$, for all $t\in\R$. But assuming so, $V=\partial_t\tilde{U}$ would 
solve
\begin{equation*}
	\begin{cases}
	\partial_t V = \partial_{xx}V-c\partial_xV & x>0,t\in\R\\
	V(0,t)=\partial_xV(0,t)=0 & x=0,t\in\R,
	\end{cases}
\end{equation*}
and therefore $V\equiv0$, so $\tilde{U}$ would be a stationary solution with range in $(0,y_1)$,
which is impossible. Hence, we can assume that (\ref{eq:48}) holds true in the case (C3).

Let $\displaystyle y_0:=\frac{1}{2}(\varphi_1(0)+\varphi_2(0))$. For all $\beta<0$, define
\begin{equation}\label{eq:50}
\psi_\beta:=\frac{\beta}{c}\left( \rme^{cx}-1\right)+y_0.
\end{equation}
Notice that $\psi_\beta$ is stationary for the evolution equation without boundary condition.
We claim that there are two functions $\varphi_3,\varphi_4\in\alpha(U)$ with 
\begin{equation*}
	\varphi_3(0)<y_0<\varphi_4(0)
\end{equation*}
satisfying 
\begin{equation}\label{eq:51}
z_{\R^+}\left(\varphi_3-\psi_\beta\right)=N=2k_1+1,\ z_{\R^+}\left(\varphi_4-\psi_\beta\right)=M=2k_2,
\end{equation}
for $\beta=g(y_0)$, and such that all zeros are simple. This is a direct consequence of (\ref{eq:48}), the 
invariance of $\alpha(U)$ under the evolution problem, and Remark \ref{convzero}. Notice that if $M=0$, 
then $\max(\psi_\beta,0)<\varphi_3\in\alpha(U)$, which would imply that $U\equiv y_1$ using the same argument given in case 
(C2), so we can assume $M\geq2$.

We differentiate between the two situations $N<M$ and $M<N$. Let us first assume that $N>M$.
Then, there exists $\beta_1>g(y_0)$, and such that (\ref{eq:51}) remains true, with all zeros simple, for $\beta=\beta_1$.
By definition of $\alpha(U)$, there exists some sequence $t_n\searrow-\infty$ such that 
\begin{equation}\label{eq:52}
   U(\cdot,t_{2k-1})\underset{k\to\infty}{\longrightarrow}\varphi_3,\ \textrm{ and }
   U(\cdot,t_{2k})\underset{k\to\infty}{\longrightarrow}\varphi_4,\ \textrm{ in }L^\infty_\mathrm{loc}(\R^+).
\end{equation}
Combining (\ref{eq:51}) and (\ref{eq:52}), we get that for all large enough $k$,
\begin{equation}\label{eq:53}
	z_{\R^+}\left( U(\cdot,t_{2k-1})-\psi_{\beta_1}\right)=N,\qquad z_{\R^+}\left( U(\cdot,t_{2k})-\psi_{\beta_1}\right)=M
\end{equation}
and all zeros are simple. Therefore, a zero was created in $[t_{2k},t_{2k-1}]$. 
Due to Lemma \ref{lemzero}, creation of a zero can only occur at $x=0$. Let $t_0\in(t_{2k},t_{2k-1})$ such that 
$U(0,t_0)=y_0=\psi_{\beta_1}(0)$. By the boundedness of $U$, the functions $\partial_tU$ and $\partial_{xx}U$ are 
uniformly bounded. Hence, since $\psi_{\beta_1}'(0)=\beta_1>\partial_xU(0,t_0)$, there exists $\eta>0$ such that 
$\displaystyle \partial_x\left( U-\psi_{\beta_1}\right)<-\eta$ on $[0,\eta]\times[t_0-\eta,t_0+\eta]$. Combining this with Lemma \ref{lemzero}
and Corollary \ref{zeroIFT}, three situations are possible:
\begin{enumerate}
	\item $U(0,t)-y_0$ has constant sign on $(s,T)$ with $s<t_0<T$, and $t\mapsto z_{\R^+}\left( U(\cdot,t)-\psi_{\beta_1}\right)$
	is non-increasing on $(s,T)$.
	\item $U(0,t)-y_0>0$ on $(s,t_0)$ and $U(0,t)-y_0<0$ on $(t_0,T)$, and $z_{\R^+}\left( U(\cdot,T)-\psi_{\beta_1}\right)<z_{\R^+}\left( U(\cdot,s)-\psi_{\beta_1}\right). $
	\item $U(0,t)-y_0<0$ on $(s,t_0)$ and $U(0,t)-y_0>0$ on $(t_0,T)$, and $\displaystyle z_{\R^+}\left( U(\cdot,T)-\psi_{\beta_1}\right)\leq z_{\R^+}\left( U(\cdot,s)-\psi_{\beta_1}\right)+1. $
\end{enumerate}
Therefore, for $k$ large enough, $z_{\R^+}\left( U(\cdot,t_{2k})-\psi_{\beta_1}\right)\leq z_{\R^+}\left( U(\cdot,t_{2k-1})-\psi_{\beta_1}\right)$, which contradicts (\ref{eq:53}). So $N>M$ is impossible.

The case $M<N$ is similar, but taking $\beta_2<g(y_0)$ instead. This concludes the proof of the first part of Proposition \ref{prop:18} in the case $c>0$ for the $\alpha$-limit set.

Concerning the $\omega$-limit set: let us assume that $U\not\equiv0$, and fix $y_0=U(0,0)\in(0,y_1)$. Define $\psi_\beta$ through (\ref{eq:50}) for $\beta<0$. Since $\alpha(U)=0$, for all $\beta<0$, there 
exists $t_\beta\ll-1$ such that
\begin{equation}\label{eq:53-1}
 z_{\R^+}\left( U(\cdot,t_\beta)-\psi_\beta\right)=1.
\end{equation}
Let $t_0=\inf\{t\leq0: U(0,t_0)=y_0\}$. Define $\beta_0=g(y_0)$, and without loss of generality, 
we can assume $t_{\beta}<t_0$ for all $\beta$. We claim that 
\begin{equation}\label{eq:53-2}
 U(\cdot,t_0)\geq\psi_{\beta_0}\ \textrm{ on }\R^+.
\end{equation}
Indeed, assume by contradiction that for some $x_0>0$, we have $U(x_0,t_0)<\psi_{\beta_0}(x_0)$. Then there exists $\beta_1=\beta_0-\delta$ with $\delta$ small enough such that 
$\displaystyle z_{\R^+}\left( U(\cdot,t_0)-\psi_{\beta_1}\right)\geq2$, and there is creation of a zero in $[t_{\beta_1},t_0]$. 
Using the same argument as above, this leads to a contradiction, and (\ref{eq:53-2}). But then, $\underline{u_0}:=\max(\psi_{\beta_0},0)$ is a sub-solution, below $U(\cdot,0)$.
Using similar arguments as in the case (C2) above, if $\underline{u}(\cdot,t)$ is the corresponding solution, then $U(\cdot,t)>\underline{u}(\cdot,t)\to y_1$, and necessarily $\omega(U)=\{y_1\}$.

The second part of Proposition \ref{prop:18} is similar to the above proof, 
with obvious modifications to rule out the possibility of elements in $\omega(u)$ between $0$ and $y_1$. 
Finally, the case $c=0$ is similar, adapting the stationary solutions $\psi_\beta(x)=\beta x+y_0$.
\end{Proof}

\begin{Remark}\label{remark:19}
 We note that the second part of Proposition \ref{prop:18} is incomplete, and the 
 natural result would be to rule out the possibility of converging to $0$. This would require 
 further analysis beyond the scope of this paper. Similar methods yield the result directly  in the non-degenerate case $g'(0)<0$. In this case, if $0\leq u_0\leq y_1$, then either $u_0\equiv0$ or 
 $\omega(u)=\{y_1\}$. In a direct adaptation of the case (C2) above we conclude that  if $u_0\not\equiv0$, then $u(\cdot,t)>0$ for any $t>0$. Then, for some $y_0$ small enough, $u(\cdot,1)<\psi_{g(y_0)}$ and 
 we can use it as a sub-solution.
\end{Remark}

\paragraph{Uniqueness of the entire solution.}
 
 In this paragraph we complete the proof of Theorem \ref{thm3} by proving that any entire 
 solution of (\ref{pb1}) with range in $(0,y_1)$ is a time translation of $U_\infty$.
The proof relies on  center manifold theory as laid out in \cite{Latushkin_DCDS08}
 for nonlinear boundary value problems. The absence of spectral gaps in the case $c=0$ makes a direct adaptation of these results, valid for $c>0$, to the case $c=0$ impossible. 
 
 We focus on the situation where the boundary term satisfies (\ref{eq:41}).
  With the transformation (\ref{symmetrization}), (\ref{pb1}) is equivalent to 
 \begin{equation}\label{eq:54}
  \begin{cases}
 \partial_t\tilde{u}+A\tilde{u}=0, & x>0,t\in\R \\
 \partial_x\tilde{u}=B(\tilde{u}) & x=0,t\in\R,
  \end{cases}
 \end{equation}
where the operator $\displaystyle A=-\partial_{xx}+\frac{c^2}{4}$ is a linear operator with domain $H^2(\R^+)$ and $B(\tilde{u})=g(\tilde{u})-\frac{c}{2}\tilde{u}$ is the nonlinear boundary condition.
The rest state $u^*\equiv0$ is an equilibrium for (\ref{eq:54}). We consider the linear system around $u^*$.
Let $A_0$ be the operator with domain $\displaystyle\mathcal{D}(A_0)=\{\varphi\in H^2(\R^+):\varphi'(0)=\left( g'(0)-\frac{c}{2}\right)\varphi(0)\}$ and $A_0=A$ on $\mathcal{D}(A_0)$.
Then, the spectrum of $A_0$ is real and given by
\begin{equation}\label{eq:spectrum}
 \sigma(-A_0)=\left(-\infty,-\frac{c^2}{4}\right]\cup\{g'(0)^2-g'(0)c\} = \sigma_\mathrm{s}\cup\sigma_\mathrm{c}.
\end{equation}
The eigenspace $E_\mathrm{c}$ is given by $\ker (A_0-\sigma_\mathrm{c}I)=\rme^{-\lambda_\mathrm{c}x}\R$ with $\lambda_\mathrm{c}=\sqrt{c^2/4+\sigma_\mathrm{c}}$ and is one-dimensional.
Let $P_\mathrm{c}$ be the spectral projection onto $E_\mathrm{c}$, $P_\mathrm{s}$ the spectral projection onto $\sigma_\mathrm{s}$. For fixed $p>3$, let us denote 
$X_0:=L^p(\R^+)$ and $X_p=W^{2(1-1/p),2}$. Applying \cite{Latushkin_DCDS08}, Theorem 5.2 in the case $g'(0)=0$ or \cite{Latushkin_JEE06}, Theorem 17 in the case $g'(0)<0$, 
there exist a function $\displaystyle \Phi\in C^1\left( P_\mathrm{c}X_0,P_\mathrm{s}X_p\right)$ with $\Phi(0)=0,\Phi'(0)=0$,
and $\varepsilon>0$ such that if $\tilde{u}(x,t)$ solves (\ref{eq:54}) with $\displaystyle \left\Vert \tilde{u}(\cdot,t)\right\Vert_{X_p}\leq\varepsilon$ for all $t<0$, then 
\begin{equation}\label{def:mcu}
 u^*+\tilde{u}(\cdot,0)\in\mathcal{M}_\mathrm{cu}:=\left\{ u^*+z_0+\Phi(z_0):z_0\in P_\mathrm{c}X_0\right\}.
\end{equation}
Notice that due to parabolic regularity, any entire solution $U(x,t)$ with $0\leq U\leq y_1$ satisfies
\begin{equation*}
 x\mapsto U(x,t)\rme^{\frac{c}{2}x}\in W^{2,q}(\R^+),\ \textrm{ for all }q\geq1.
\end{equation*}
By definition of $\Phi$, since $E_\mathrm{c}$ is one-dimensional, $\mathcal{M}_\mathrm{cu}$ is a one-dimensional center unstable manifold.
For a fixed $y_0$, consider the entire solution defined in (\ref{eq:40}). Then $U_0(\cdot,t)$ converges to $0$ by the above as $t\to -\infty$. Therefore, 
up to a smaller $\varepsilon$ if necessary, there exists some $T\in\R$ such that
\begin{equation}\label{eq:56}
 \mathcal{M}_\mathrm{cu}\cup B_{X_p}(u^*,\varepsilon)\cup\{\tilde{u}:\tilde{u}\geq u^*\}=\left\{ U_\infty(\cdot,t)\rme^{-\frac{c}{2}\cdot}:t<T\right\}.
\end{equation}
Let $U(x,t)$ be any other entire solution of with $0<U<y_1$. By Proposition \ref{prop:18}, $U(\cdot,t)\to0$ as $t\to-\infty$ in $C^2_\mathrm{loc}(\R^+)$.
Therefore, there exists $T_1\in\R$ such that $V(x,t)=U(x,t)\rme^{\frac{c}{2}x}$ satisfies $\displaystyle \left\Vert V(\cdot,t)-u^*\right\Vert_{X_p}\leq\varepsilon$, for all $t<T_1$, 
and so by (\ref{def:mcu}-\ref{eq:56}) necessarily $U(\cdot,t)=U_\infty(\cdot,t-T_2)$ for some $T_2\in\R$, for all $t<T_1$. Consequently, any entire solution with range in $(0,y_1)$
is a time-translation of $U_\infty$. This concludes the proof of Theorem \ref{thm3}. 

\section{Purely diffusive dynamics}\label{s:5}

 
We conclude the proofs of the main results stated in the introduction. 
Having established Proposition \ref{p:het0} in the previous section, we only need 
to prove Proposition \ref{prop:4}. We also present several results and observations 
on possible uniqueness statements, striving to describe in more detail asymptotics of heteroclinic solutions. 


\paragraph{Asymptotic behavior in the case $g>0$ --- proof of Proposition \ref{prop:4}.}
Assuming $0<\gamma_1\leq g\leq \gamma_2$, we proceed via a direct computation. Let $u_0$ be any initial condition in $L^\infty(\R^+)$.
Then, the comparison principle given in Lemma \ref{comparison} states that
\begin{equation}\label{eq:70}
 u_2(\cdot,t)<u(\cdot,t)<u_1(\cdot,t),\ t>0
\end{equation}
where $u_{1,2}$ are solutions of 
\begin{equation}\label{eq:71}
  \begin{cases}
  \partial_t u_i = \partial_{xx}u_i, \qquad & x>0,\ t>0\\
  \partial_x u_i = \gamma_i, \qquad & x=0,\ t>0, \\
  u_i(x,0)=(-1)^{i-1}\left\Vert u_0\right\Vert_\infty,\qquad & x\geq0.
 \end{cases}
\end{equation}
The solutions of (\ref{eq:71}) are explicit and given by
\begin{equation}
 u_i(x,t)=(-1)^{i-1}\left\Vert u_0\right\Vert_\infty-2\gamma_i\int_0^t\frac{\rme^{-\frac{x^2}{4(t-s)}}}{\sqrt{4\pi(t-s)}}\rmd s,
\end{equation}
which immediately gives Proposition \ref{prop:4}.

\paragraph{Dynamics between stationary solutions.}
Unlike Proposition \ref{prop:17} and \ref{prop:18}, the uniqueness of the entire solution 
in Theorem \ref{thm3} does not extend to the case $c=0$. There are two main obstacles to this.
First, the a priori estimates in Proposition \ref{prop:18} are sufficient in the case $c>0$
since we work in exponentially weighted spaces; the same information does not appear to be sufficient 
in the case $c=0$ to control the behavior at $x=+\infty$. Second, exponential weights provide a spectral gap in the case $c>0$, whereas for $c=0$ we have  continuous
spectrum up to the origin, no point spectrum when $g'(0)\geq 0$. We can nevertheless construct unique entire solutions 
with prescribed convergence rates as $t\to-\infty$, using either strong unstable manifold theory or scaling variables \cite{Wayne_ARMA97}, to overcome continuous spectrum at the origin. 
Throughout we assume  (\ref{eq:41}), and keep the notation of the previous section with the formalism associated with (\ref{eq:54}).

\subparagraph{Linear behavior near the origin.}
Let us first assume that the local behavior around $0$ is 
given by 
\begin{equation}\label{eq:74}
 g(y)=-\gamma y+\rmO(y^p)\textrm{ for some }p>1.
\end{equation}
The asymptotics of solutions for $t\to -\infty$ was 
described in Proposition \ref{prop:18} and Remark \ref{remark:19}. We shall discuss uniqueness in classes of exponentially decaying functions, here. 
The relevant  linear operator is $A_0=-\partial_{xx}$ with domain $\displaystyle\mathcal{D}(A_0)=\{\varphi\in H^2(\R^+):\varphi'(0)=-\gamma\varphi(0)\}$. It admits a spectral gap:
\begin{equation}\label{eq:spectrum2}
 \sigma(-A_0)=\left(-\infty,0\right]\cup\{\gamma^2\} = \sigma_\mathrm{s}\cup\sigma_\mathrm{u}.
\end{equation}
The eigenspace associated to the positive eigenvalue is one-dimensional, given by 
$\displaystyle E_\mathrm{u}=\rme^{-\gamma x}\R$. Considering the functions spaces $X_0, X_p$ defined in the previous 
section, and letting $P_\mathrm{s},P_\mathrm{u}$ be the spectral projection onto $\sigma_\mathrm{s},\sigma_\mathrm{u}$, there exists 
$\Phi\in C^1(P_\mathrm{u}X_0,P_\mathrm{s}X_p)$ with $\Phi(0)=0,\Phi'(0)=0$ such that the unstable manifold 
\begin{equation*}
 \mathcal{M}_\mathrm{u}=\left\{z_0+\Phi(z_0):z_0\in P_\mathrm{u}X_0\right\}
\end{equation*}
is invariant under the evolution problem as long as the solution remains small. Moreover, 
due to the spectral gap (\ref{eq:spectrum2}) for 
all $\delta\in(0,\gamma^2)$, if $U(x,t)$ is a solution coming from $0$ as $t\to-\infty$,
then for some $T\ll-1$,
\begin{equation}\label{eq:76}
 \underset{t\leq0}{\sup}\left\{\left\Vert U(\cdot,t)\rme^{-(\gamma^2-\delta)t}\right\Vert_{X_p}<\infty\right\}\Longleftrightarrow\left( U(\cdot,t)\in\mathcal{M}_\mathrm{u},\textrm{ for }t<T\right).
\end{equation}
Since $P_\mathrm{u}X_0=E_\mathrm{u}$ is known and one-dimensional, the asymptotics (\ref{eq:76})
is satisfied by a one parameter family of solutions, of the form
\begin{equation*}
 U(x,t)=\rme^{-\gamma x+\gamma^2(t+T)}+\tilde{\Phi}(T)
\end{equation*}
where $\tilde{\Phi(T)}=\Phi(\rme^{\gamma^2T-\gamma\cdot})$.

\subparagraph{Quadratic behavior around the origin.}
In this critical case, we assume that the behavior of $g$ near the origin is given by 
\begin{equation}
 g(y)=-y^2+g_1(y),\ \textrm{ with }g_1(y)=\rmO(y^3).
\end{equation}
The spectrum of the linearized operator is simply $\R^-$, and we resort to similarity variables in order to understand the algebraic rates of decay. 
Since we are interested in asymptotics for $t\to-\infty$, we assume that $t<-1$. 
For an ancient solution $U$ of (\ref{pb:ceq0}), consider the change of variables 
\begin{equation}\label{eq:77}
 \tau=-\log(-t),\ \xi=\frac{x}{\sqrt{-t}},\ \textrm{ and \ }U(x,t)=\frac{1}{\sqrt{-t}}V(\xi,\tau).
\end{equation}
Then the resulting equation for $V$ reads, using
$\displaystyle \eta=1/\sqrt{-t}:=\rme^{\tau/2}$,
\begin{equation}\label{eq:79}
 \begin{cases}
\partial_\tau V=\partial_{\xi\xi}V-\frac{\xi}{2}\partial_\xi V-\frac{V}{2}=-LV, & \xi>0,\ \tau<-1 \\
\partial_\xi V=-V^2+\rmO(\eta V^3), & \xi=0,\ \tau<-1 \\
\partial_\tau\eta=\frac{\eta}{2}.
 \end{cases}
\end{equation}
The stationary equation can be solved explicitly in $\xi>0$, which results in an algebraic equation at $\xi=0$ that we solved, finding a unique equilibrium
$(V,\eta)=(V^*,0)$ with 
\begin{equation}\label{eq:80}
 V^*(\xi)=\frac{1}{\sqrt{\pi}}\rme^{\frac{\xi^2}{4}}\textrm{erfc}(\xi),
\end{equation}
in the ambient space $H^2(\R^+)\oplus\R$. Linearizing at this stationary 
solution, the operator becomes $L_0=L$ on $\mathcal{D}(L_0)=\{\varphi:\varphi'(0)=-2V^*(0)\varphi(0)\}$. In spaces of sufficiently localized functions, 
the spectrum of this operator is discrete, as one can quickly infer from the methods in \cite{Wayne_ARMA97}. We found that 
\begin{equation*}
 \sigma(-L_0)=\sigma_\mathrm{s}\cup\sigma_\mathrm{u}
\end{equation*}
with $\sigma_\mathrm{s}\subset(-\infty,\lambda_\mathrm{s}]$, $\lambda_\mathrm{s}=-1.23162\ldots$, and $\sigma_\mathrm{u}=\{1\}$. Notice also that the linearization 
of (\ref{eq:79}) admits another positive eigenvalue in the $\eta$ direction, for which the dynamics 
is trivial. Let $\varphi_1$ be the eigenfunction associated to $\sigma_\mathrm{u}$. Using again the methods in \cite{Wayne_ARMA97}, 
there exists $\Phi\in C^1(\R^2,P_\mathrm{s}X_p)$ with $\Phi(0)=0$ and $\Phi'(0)=0$, and $\varepsilon>0$, such that the unstable manifold 
\begin{equation}
 \mathcal{M}_\mathrm{u}=\left\{ \left( V^*+A\varphi_1+\Phi(A,\eta),\eta\right):(A,\eta)\in(-\varepsilon,\varepsilon)\times[0,\varepsilon)\right\}
\end{equation}
is invariant for the evolution problem (\ref{eq:79}) for $\tau<-1$. Going back to the $(x,t)$
variables, we get the existence of an entire solution 
\begin{equation*}
 U(x,t)=\frac{1}{\sqrt{-t}}V^*\left( x/\sqrt{-t}\right)+\tilde{\Phi}(t)
\end{equation*}
where $\displaystyle \tilde{\Phi}(t)=\frac{1}{\sqrt{-t}}\Phi(0,1/\sqrt{-t})=O\left((-t)^{-3/2}\right).$ 
This entire solution is unique within the class of solutions that are small and bounded as $t\to -\infty$ in the scaling variables.

\section{Discussion}\label{s:d}

Summarizing, we established existence, and, within reasonable classes of functions, uniqueness of entire solutions, thus characterizing the asymptotic behavior  in this class of advection-diffusion equations with nonlinear flux and gauge symmetry. Our results are schematically summarized in Table \ref{t:1}. 
\begin{table}
\begin{tabular}{>{\centering\arraybackslash}p{1.2in}>{\centering\arraybackslash}p{1.2in}>{\centering\arraybackslash}p{1.2in}>{\centering\arraybackslash}p{1.2in}}
 \multicolumn{4}{c}{\textbf{positive growth rate} $\mathbf{c>0}$} \\[0.1in]
\hline\\[-0.1in]
& $g>0$ & $g\geq 0$ & $g\gtrless 0$\\[0.1in]
\hline\\[-0.1in]
\raisebox{0.2in}{dynamics} &\includegraphics[width=.5in]{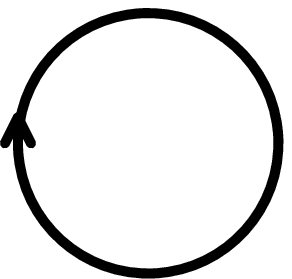} &\includegraphics[width=.5in]{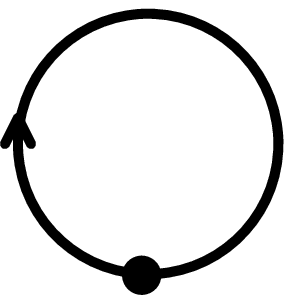} &\includegraphics[width=.5in]{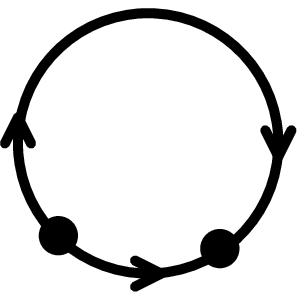} \\
\hline\\[-0.1in]
existence & \checkmark \ \ Thm~\ref{thm1}&\checkmark\ \  Thm~\ref{prop2}&\checkmark\ \  Thm~\ref{thm3}\\[0.1in]
\hline\\[-0.1in]
uniqueness  & \checkmark \ \ Thm \ref{thm1}&\checkmark\ \  Thm~\ref{prop2}&\checkmark\ \  Thm~\ref{thm3}\\[0.1in]
\hline\\[-0.1in]
stability & \checkmark \ \ Thm \ref{thm1}&\checkmark\ \  Thm~\ref{prop2}&\checkmark\ \  Thm~\ref{thm3}\\[0.1in]
\hline\\[-0.1in]
 interpretation & resonant  growth & \multicolumn{2}{>{\centering\arraybackslash}p{2.4in}}{compatible zero strain selection,\newline phase slips  from small perturbations}\\[0.1in]
\hline\\[-0.1in]
\end{tabular}\\[0.4in]

\begin{tabular}{>{\centering\arraybackslash}p{1.2in}>{\centering\arraybackslash}p{1.2in}>{\centering\arraybackslash}p{1.2in}>{\centering\arraybackslash}p{1.2in}}
 \multicolumn{4}{c}{\textbf{bounded solutions}, \textbf{zero growth rate} $\mathbf{c=0}$} \\[0.1in]
\hline\\[-0.1in]
& $g>0$ & $g\geq 0$ & $g\gtrless 0$\\[0.1in]
\hline\\[-0.1in]
\raisebox{0.2in}{dynamics} &\raisebox{0.08in}{\includegraphics[width=.5in]{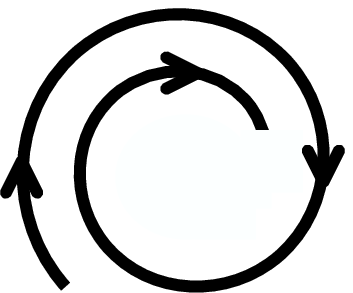}} &\includegraphics[width=.5in]{sn.eps} &\includegraphics[width=.5in]{het.eps} \\
\hline\\[-0.1in]
existence & \checkmark \ \ Prop \ref{prop:4}&\checkmark \ \ Prop \ref{p:het0}&\checkmark \ \ Prop \ref{p:het0}\\[0.1in]
\hline\\[-0.1in]
uniqueness  & $\times$ & (\checkmark) \ \ Prop \ref{p:het0} &(\checkmark) \ \ Prop \ref{p:het0}\\[0.1in]
\hline\\[-0.1in]
stability & $\times$ &(\checkmark) \ \ Prop \ref{p:het0} &(\checkmark)  \ \ Prop \ref{p:het0}\\[0.1in]
\hline\\[-0.1in]
 interpretation &  diffusion $\sim\sqrt{t}$ toward compatible strain & \multicolumn{2}{>{\centering\arraybackslash}p{2.4in}}{compatible zero strain selection, \newline diffusive phase slips from small perturbations}\\[0.1in]
\hline\\[-0.1in]
\end{tabular}
\caption{Schematic representation of our results phase portraits up to gauge symmetry. See the main theorems for precise statements, in particular concerning the (partial) uniqueness results in the case $c=0$.}
\label{t:1}
\end{table}
While the table shows apparent strong similarities between the cases $c>0$ and $c=0$, there are important differences. First, there is in fact a family of tables for the case $c=0$, parameterized  by $k$,  the asymptotic strain $u\sim kx$, $x\to \infty$. Second, rates are much slower, diffusive, rather than exponential. 

\paragraph{Back to crystal growth.}
Interpreting Table \ref{t:1} in light of the motivation from crystal growth, we think of the slope $u_x$ as the wavenumber or strain and the value $u$ as the phase. The trivial case, origin of our homotopies, where $g\equiv g_0$ is constant corresponds to boundary conditions that are compatible with one particular strain, only, for arbitrary phases. From this perspective, our results account for non-adiabatic effects, where the effect of boundaries depends on the phase at the boundary, that is, the crystalline microstructure at the boundary cannot be averaged. In the case $g>0$, we establish that growth processes will be \emph{resonant},  periodic with minimal period. This in particular excludes subharmonic bifurcations, when solutions emerge who are periodic up to multiples $2\ell\pi$, $\ell>1$, of the gauge symmetry $2\pi$.

In this context, the interesting quantity is the relation between the speed $c$ and the crystal strain $k=\lim_{x\to\infty}u_x$, which is related to 
the period through the simple relation $ck=\omega=2\pi/T$. The coarse nature of our results does not provide a quantitative description 
of this speed-strain relation, but rather establishes existence, uniqueness, and smoothness of a relation $k(c)$. We comment below on more quantitative results and contexts where more complicated, including subharmonic, dynamics can arise. 

The case when $g$ possesses zeros is not immediately relevant for the case of non-trivial patterns but does arise in systems with trivial constant solutions such as  the complex Ginzburg-Landau equation \cite{beekie,ms}. Zeros of $g$ correspond to constant phase solutions which are naturally ``compatible'' with domain growth, such that the convergence to these states in the case $c>0$ can be interpreted as a selection of compatible, zero strain. In the stationary case, $c=0$, the asymptotic strain is conserved. Imposing then, for instance, zero strain, the system either relaxes to a phase with zero strain, $g(u_0)=0$, or, when $g$ does not possess zeros, diffusively approaches a state with (extremal) strain; see Proposition \ref{prop:4}.

\paragraph{Open questions.}
There are clearly many more subtle questions one could ask, pertaining for instance to global attractivity in the case $g>0$, 
asymptotics in the case $g>0$, $c=0$, or uniqueness in larger classes of functions in the case $c=0$, $g$ sign-changing.  
A question of more general interest arises when describing ancient solutions in the case $g(u_0)=0$, $g'(u_0)<0$. Besides the exponentially 
decaying solutions that we analyzed here, one would like to exclude solutions that decay inside of a center-stable manifold. 
It is however not clear how to use the equivalent of  scaling variables after projecting onto this codimension-one manifold. 

\paragraph{Connecting the entire solutions across parameters.} Our picture may appear somewhat fractured, and we are indeed 
missing a more cohesive description of how entire solutions limit onto each other as both $c$ and $\vartheta$ are varied. 
Fixing $c>0$, we have a fairly complete picture, analogous to the classical saddle-node-on-limit-cycle (SNIC) bifurcation in 
dynamical systems, where periodic orbits limit on a saddle-node homoclinic, which then breaks up into two heteroclinic 
connections forming an invariant circle. Our picture is slightly more complicated due to the presence of essential spectrum, 
which leads to complications such as the rather weak convergence in $L^\infty_\mathrm{loc}$. A more refined description would 
include the release of an error-function type kink as heteroclinic solutions converge for $t\to+\infty$. Information clearly 
is more limited, yet, in the case $c=0$, and we have not attempted to attach the less selective dynamics that occur for $c<0$, 
when initial conditions are advected towards the boundary and patterns are annihilated. In this case, one would eliminate 
essential spectrum by using exponentially growing weights, thus preserving the asymptotic growth and preventing the release 
of error-function type kinks. We expect periodic solutions with prescribed strain $k$, rather than the selected strain $k(c)$ 
in the case $c>0$. 

\paragraph{Higher space dimensions.} Many of our techniques here can be adapted to the multi-dimensional setting, 
studying advection-diffusion in $x_1>0,y\in\R^n$, with transport speed $c\in \R^n$ away from the boundary, $c_1\geq 0$. 
One can then prescribe a lateral strain $k_y$ through a periodicity assumption and mimic for instance the construction 
of relative periodic orbits, here. We intend to pursue these questions in forthcoming work. 

\paragraph{Quantitative speed-strain relations.} Quantitative results were obtained in a formal setting in \cite{beekie}. 
The regions of interest there were $0<c\ll 1$ and $c\gg 1$. The former case leads to a rather involved singular 
perturbation problem and asymptotics $k\sim \mathrm{min}\,g -\sqrt{2}\zeta(1/2) \sqrt{c}+\rmO(c^{3/4})$, 
where $\zeta$ is the Riemann $\zeta$-function. 
The latter case, $c=\frac{1}{\varepsilon^2}$, can be thought of as a large-advection limit after scaling  $c t=\tau$,
\[
 \partial_\tau u = \varepsilon^2\partial_{xx}-u_x,
\]
with limit $u_\tau=u_x$ in $x\geq 0$, hence $u_t=-g(u)$  at $x=0$. This can be made rigorous expanding 
the Dirichlet-to-Neumann operator with periodic boundary conditions, or, equivalently, constructing slow manifolds 
in spatial dynamics \cite{scheel94}. At $\varepsilon=0$, infinite speed, we can solve the differential equation at 
the boundary explicitly to find the period $T=\int_0^{2\pi}g(u)^{-1}\rmd u$, and therefore $\omega=\left(\dashint g^{-1}\right)^{-1},$ 
the harmonic average of the boundary flux. Expanding the Dirichlet-to-Neumann operator further, one can also compute higher-order corrections; see \cite{beekie}.

\paragraph{Transport versus memory.} The simple form of our equation in $x>0$ would allow us to reduce the dynamics to a 
pseudo-differential equation on the boundary, introducing in particular infinite memory into the dynamics. 
In the purely diffusive case, this reduction has been exploited in the literature, giving rise to the field of fractional differential equations. 
The memory kernel is increasingly localized as the speed increases, which can intuitively be understood as reducing 
the effect of the pattern away from the boundary by increasing the advection speed. In the limit of infinite speed, 
the equation then becomes a local-in-time differential equation. The corrections in $k(c)$ for finite $\varepsilon=1/\sqrt{c}$ 
can then be understood as incorporating moment approximations of this memory kernel; see for instance \cite{fayescheel} for such 
reductions and expansions for nonlocal kernels. The SNIC bifurcation, at $|\vartheta|=1$,  would then be understood in 
the context of saddle-node bifurcation in the presence of memory, or, in the extreme case $c=0$ as a saddle-node bifurcation in a fractional differential equation.

\paragraph{Strain-speed relations beyond phase-diffusion.}
Turning to the original motivation from striped phases, one would like to know how much of the present results or techniques 
transfer to systems such as the Swift-Hohenberg equation. Technically, the methods here are based on comparison principles
and therefore not applicable. One can however envision scenarios such as in \cite{goh3,weinburd} where phase-diffusion equation 
with effective boundary could be derived rigorously and the results obtained here would then indeed imply resonant crystal growth 
in the Swift-Hohenberg equation. The results in \cite{beekie} demonstrate numerically that working on the level of the phase-diffusion 
equation while extracting the nonlinear flux $g$ quantitatively in the form of a strain-displacement relation can give accurate 
quantitative predictions for a large class of pattern-forming equations.

\paragraph{Phenomena beyond phase-diffusion.} 
On the other hand, there are numerous phenomena, just within the Swift-Hohenberg equation, that cannot be analyzed in this simple scalar context, including detachment and \cite{goh1,goh2} and  wrinkling \cite{avery}. In fact, more complex phenomena arise precisely when the phase-diffusion approximation fails to capture the growth dynamics, due to the relevance of amplitude variations or other instabilities. The simplest example are side-band instabilities, in which case the sign of the diffusion constant changes and fourth-order differential operators are needed to stabilize dynamics. Such side-band instabilities are a key ingredient in two-dimensional growth scenarios when the orientation of stripes is perpendicular to the boundary. It was demonstrated in \cite{avery} that the growth dynamics are well approximated by a Cahn-Hilliard equation for the phase gradient together with an advection term and effective boundary condition at the growth interface. Different from the scenario here, the simplest ``harmonic'' or ``resonant'' growth destabilizes in a saddle-node bifurcation on a limit-cycle, giving rise to more complex patterns, corresponding to crystals with superstructure, specifically wrinkled stripes. Translated into the context of our results, the Cahn-Hilliard setting forces solutions to possess asymptotically periodic strain, with $\Phi_x$ alternating periodically between values close to $\pm 1$ rather than approaching a finite limit. The interfaces between regions of different strain, known as kinks, interfaces, or defects in different contexts, are stationary in a steady frame and can be interpreted as building blocks for superstructures in crystals. Defect nucleation during crystal growth has also been documented in the example of the complex Ginzburg-Landau equation in \cite{beekie}, where amplitude defects emerged for moderate growth speeds which pushed the selected strain into an Eckhaus-unstable regime.

\bibliographystyle{abbrv}
\footnotesize
\bibliography{pdg}

\Addresses

\end{document}